\documentclass[a4paper,12pt]{article}
\usepackage{a4wide,psfig}

\title
{A route to computational chaos revisited:
noninvertibility and the breakup of an invariant circle}

\author
{Christos E. Frouzakis  \\
   Combustion Research Laboratory,
   Paul Scherrer Institute\\
   CH-5232, Villigen, Switzerland
\and
Ioannis G. Kevrekidis 
\footnote{Corresponding author, Department of Chemical Engineering,
Princeton University, Princeton, NJ 08544, Phone: (609) 258 2818,
Fax: (609) 258 0211, e-mail: yannis@arnold.princeton.edu}\\
  Department of Chemical Engineering\\
and Program in Applied and Computational Mathematics,\\
  Princeton University, Princeton, NJ 08544
\and
Bruce B. Peckham\\
Department of Mathematics and Statistics\\
University of Minnesota Duluth, Duluth, MN 55812}

\begin{document}
\vspace{-.9in}
\begin{small}
\maketitle
\end{small}
\vspace{-.3in}
\begin{abstract}
In a one-parameter study of a noninvertible family of maps of the plane
arising in the context of a numerical integration scheme, Lorenz
studied a sequence of transitions from an attracting fixed point to 
``computational chaos.''
As part of the transition sequence, he
proposed the following as a possible
scenario for the breakup of an invariant circle:
%
the invariant circle develops regions of increasingly sharper curvature
until at a critical parameter value it develops cusps;
beyond this parameter value, the invariant circle fails to persist,
and the system exhibits chaotic behavior on an invariant set with loops \cite{Lo2}.
We investigate this problem in more detail and show that the invariant
circle is actually
destroyed in a global bifurcation before it has a chance to develop cusps.
Instead, the
%
global unstable manifolds of saddle-type periodic points
are the objects which develop cusps and subsequently ``loops'' or ``antennae.''
The one-parameter study is better understood when embedded in
the full two-parameter space and viewed in the context of
the two-parameter Arnold horn structure.
%
Certain elements of the interplay of noninvertibility 
with this structure, the associated invariant circles,
periodic points and global bifurcations are examined.
\end{abstract}


\eject

\section{Introduction}\label{s-intro}

In an insightful 1989 paper entitled ``Computational Chaos: a prelude to
computational instability" E. N. Lorenz \cite{Lo2}
reported on a one-parameter
computational study of the dynamics of the {\it noninvertible} map
\begin{eqnarray}
L_{(a,\tau)}:R^2 \rightarrow R^2 = \left\{ \begin{array}{l}
                      x_{n+1}=(1+a\tau)x_n-\tau x_n y_n \\
                      y_{n+1}=(1-\tau)y_n+\tau x_n^2, 
      \end{array} \right.   \label{map}
\end{eqnarray}
which arises from a simple forward Euler integration scheme 
($\tau$ being the time step of the integration) of the 
two coupled nonlinear ODEs
\begin{eqnarray}
     dx/dt &=& ax - xy \nonumber \\ 
     dy/dt &=& -y +x^2.  \label{ode}
\end{eqnarray}
These ODEs are obtained by starting with the familiar 
Lorenz system \cite{Lo1}:
$\dot{x}=-\sigma (x-y), \dot{y}=-xz+\rho x -y, \dot{z}=xy-\beta z$,
letting $\sigma \rightarrow \infty$, and rescaling the variables and 
the remaining parameters.
Lorenz fixed the value of $a$ at $0.36$, and varied the time step $\tau$.
Although the corresponding differential equations (eq.~(\ref{ode})) exhibit
an attracting equilibrium point, computer simulations indicated that
the discrete approximation (eq.~(\ref{map}))
progressed from exhibiting an attracting fixed point, to an attracting
invariant circle (IC), to a chaotic attractor (termed {\it computational chaos})
as the time step $\tau$ was gradually increased.
%
%
Sequences of bifurcations similar to those described by Lorenz have been
observed for other noninvertible maps of the plane \cite{MiraBook3, RAK},
suggesting
that this might be a universal ``noninvertible route to chaos.''
Some of the bifurcations along the route
--- local changes of stability, homoclinic and heteroclinic tangencies,
crises caused by an attractor interacting with its basin boundary ---
are also observed in families of {\it invertible} maps.
New bifurcations, however, unique to noninvertible families,
as well as ``invertible bifurcations with noninvertible complications,''
are also observed along the route.

In this paper we perform a more detailed numerical study of Lorenz's family,
but focus primarily on a narrow range of parameters which includes the breakup
of the invariant circle and the identification of $\tau_{CD}$ (the greatest
lower bound on the values of the parameter $\tau$ for which the corresponding
maps have ``chaotic dynamics''), and $\tau_{chaos}$,
the greatest lower bound on the values of the parameter $\tau$
for which the corresponding maps exhibit a ``chaotic attractor.''
Lorenz was interested in $\tau_{chaos}$
(labelled $\tau_b$ in \cite{Lo2}) as an indication
of how poorly the Euler map
$L_{(0.36,\tau)}$ approximated the
original
differential equation.
%
%
He observed that, in the transition from smooth
IC to chaotic attractor, the process began with
the circle developing features with increasingly high curvature.
There appeared to be a
critical $\tau$ value at which the ``IC developed cusps.''
After this value, the ``chaotic attractor'' suggested by computer simulations
appeared to have loops and a Cantor-like structure.
It certainly was no longer a topological circle.
He claimed (correctly)
that if this scenario did in fact occur, then this critical $\tau$ value
had to be an upper bound on $\tau_{chaos}$.

%
Our numerical investigations
suggest that a smooth IC does not persist all the way
to a ``cusp'' parameter.
Instead, the IC is destroyed in a heteroclinic tangency
initiating the crossing of a branch of the stable manifold
and a branch of the unstable manifold of a periodic (here period-$37$)
saddle point.
(As for invertible maps, chaotic dynamics are apparently present during this
crossing, so this tangency is an upper bound on $\tau_{CD}$.)
It is the unstable manifold that subsequently develops cusps at a critical
parameter value we label $\tau_{cusp}$, 
and loops beyond that value.
A second manifold crossing results in the apparent appearance of the chaotic
attractor, with loops inherited from the unstable manifold.
We conclude that $\tau_{cusp}$ is not a single isolated bifurcation value
separating smooth IC attractors from chaotic attractors.

Rather, $\tau_{cusp}$ is strictly above $\tau_{CD}$ and strictly below
$\tau_{chaos}$.  All three parameter values are part
of the transition from a smooth IC attractor
to a chaotic attractor.
Because a cusp on an IC with an irrational rotation number would
necessarily force a dense set of cusp points on the IC, it would make
the existence of the IC itself unlikely.
We therefore expect the existence of an entire transition interval of
parameters,
as opposed to a single bifurcation parameter, to be the generic scenario
for the transition between smooth IC and chaotic attractor in the presence of
noninvertibility.
%

Mechanisms for the breakup of ICs
are relatively well established for {\it invertible}
maps of the plane (see for example \cite{ACHM,Chenciner, Pthesis}).
In the invertible scenario which appears to us to
most closely match the noninvertible
scenario of this paper, an IC is born in a Hopf
(also called Neimark-Sacker) bifurcation,
grows in size, coexists with a periodic orbit (after a saddle-node birth
of the periodic orbit off the IC), is destroyed in a first
global manifold crossing
(also referred to in the literature as a crisis \cite{OttBook}
where the IC 
collides with its basin boundary),
and is reconstituted after a second global manifold crossing.
The periodic orbit which persists through the
destruction and reconstitution of the IC typically switches from
outside (inside) the IC before it is destroyed to inside (outside)
the IC after it is reconstituted.
Chaotic invariant sets necessarily exist only during the manifold crossings.
In contrast, the breakup of the IC in the noninvertible scenario
is part of a transition to a chaotic attractor.
In particular,
the ``chaotic attractor with loops,''
which appears in this noninvertible case, is not a feature which appears in
invertible maps.
In fact, cusps and loops on iterates of smoothly embedded curves
are not possible for smooth invertible maps, but
are common features on (invariant) curves of noninvertible maps \cite{FGKM}.
Further, bifurcations such as the manifold crossings (crises)
can have
additional complications due to the presence of noninvertibility.
For example, the unstable ``manifolds'' involved in the global
crossings may have self intersections and cusps.
Or stable and unstable ``manifolds'' may cross transversely at one 
homoclinic point,
but fail to preserve transversality at other points along the homoclinic
orbit.

In any case,
the transition mechanism for the maps in eq.~(\ref{map}) is different
from the invertible transitions,
and the noninvertible nature of the map plays an important role:
the IC --- for parameter values for which
it exists --- interacts with (actually crosses) the locus on which
the Jacobian of the linearized map becomes singular.
This locus (and its images and preimages)  is crucial in
organizing the dynamics of noninvertible maps; it is termed
``critical curve" in the pertinent literature 
\cite{MiraBook1,MiraBook2,MiraBook3,MiraBook4} and constitutes the
generalization (in two dimensions) of the critical point in
unimodal maps of the interval \cite{Feigenbaum,ColletEckmann}.
In particular, as explained in \cite{FGKM}, the angle of intersection
between a saddle unstable manifold and the critical curve underpins the
transition of the local image of the manifold
from being smooth and injective to nonsmooth and injective (with a cusp)
to smooth but not injective (with loops).

As for invertible Hopf bifurcations (even in the noninvertible setting the
dynamics near a Hopf bifurcation are necessarily locally invertible), an
understanding of Lorenz's $a=0.36$ one-parameter family is possible only
by viewing it in the two-parameter Arnold horn context
\cite{Arnold, Chenciner, ACHM, P, PFK}.
We do this with $L_{(a,\tau)}$, using the second parameter ($a$)
already provided in eq.~(\ref{map}).
This leads to an examination of the internal Arnold resonance
horn structure, to be compared to and contrasted with the invertible case.
The interiors of the
resonance horns have features which vary significantly from those in
the horns of invertible maps (see preliminary results in \cite{GFK, GAK}).
Although our understanding of this internal structure is far from
complete, examination of the partial picture still provides
insight into the bifurcations observed in the Lorenz one-parameter cut.
Developing computational tools to further investigate the internal
structure is part of our ongoing research.

In this paper the observations of the transition in \cite{Lo2} are
briefly summarized in Sec.~\ref{s-observations}, and then revisited in
detail in Sec.~\ref{s-revisit}.
In Sec.~\ref{s-larger}, we place our one-parameter cut in context by
examining the ``resonance horns'' in the larger two-parameter space.
We discuss related noninvertible issues, including computational challenges,
in Sec.~\ref{s-discussion}, and state final conclusions in 
Sec.~\ref{s-conclusions}.

\section{Observations}\label{s-observations}

\begin{figure}
\centerline{ \psfig{figure=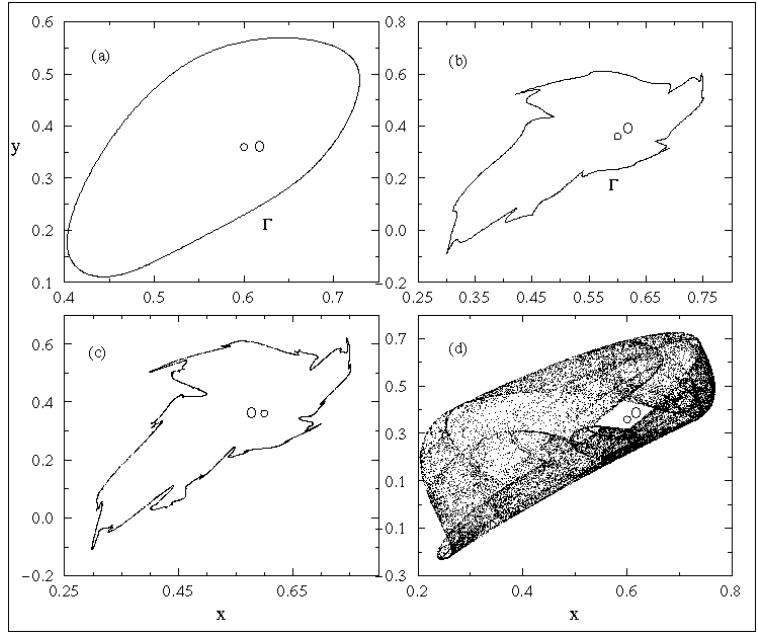,height=20pc} }
\caption{ Phase portraits of the map $L$ for $a=0.36$ and
(a) $\tau=1.55$, smooth IC 
(b) $\tau=1.775$, regions of sharp curvature appear on the IC
(c) $\tau=1.785$, self-intersections appear (d) $\tau=1.91$, ``fully chaotic''
attractor
\label{observ}}
\end{figure}

Figure \ref{observ} briefly summarizes the relevant observations in \cite{Lo2}.
For $a$ fixed at $0.36$, the differential equation
(eq.~(\ref{ode})) has
an attracting equilibrium point at $(x,y)=(0.6, 0.36)$.  Its basin of
attraction is the right half plane.  There is a symmetric attracting
equilibrium for
the left half plane at $(-0.6, 0.36)$.  We will deal only with the right half
plane attractors in this paper.
The corresponding maps $L_{(0.36,\tau)}$ all have a fixed point at
$(x,y)=(0.6, 0.36)$.
In a one-parameter cut with respect to $\tau$, the fixed point is stable
for small $\tau$,
but loses its stability via Hopf bifurcation at
$\tau = 1.38889$.
The resulting attractor  is initially a smooth IC
(fig.~\ref{observ}a, $\tau = 1.55$) but then
develops progressively sharper features (fig.~\ref{observ}b, $\tau = 1.775$,
where the circle is still smooth even though it appears to have cusps),
then self-intersects, suggesting chaos 
(fig.~\ref{observ}c, $\tau = 1.785$), and eventually exhibits clear signs
of chaotic behavior (fig.~\ref{observ}d, $\tau = 1.91$).  Lorenz
supports his claim of chaos with the computation of a positive Lyapunov
exponent for $\tau=1.91$ \cite{Lo2}.
The last two figures are 
unusual in that the attractor apparently shows self-intersections, a clear
sign of noninvertibility.
All four figures  were created by following a single orbit, after dropping
an initial transient, from an initial
condition near the fixed point.
We use the term {\it attractor} loosely to denote the object
generated by the computer via such a simulation.  
More formal definitions appear in the next section.

Lorenz proposed the following scenario as the simplest possible
transition of the attractor from a smooth IC
to a chaotic attractor.
As a parameter is increased ($\tau$ in this example) and the IC crosses
successively farther over the
critical curve (called ``$J_0$'' below), causing it to
develop successively sharper features.
Simultaneously, its rotation number changes (increasing with $\tau$ in
this case).
By restricting to $\tau$ values for which the IC is quasiperiodic
(and avoiding periodic lockings, defined below in Sec.~\ref{ss-prelims}),
one obtains
a sequence of $\tau$ values, each corresponding to a map with a smooth
quasiperiodic IC, and limiting to $\tau_{chaos}$, the lower limit on
chaotic behavior.
At $\tau_{chaos}$, the attractor, which may or
may not still be a topological circle, would develop cusps.
Beyond $\tau_{chaos}$, the attractor would be chaotic, or at least have parameter
values accumulating to $\tau_{chaos}$ from above for which the corresponding maps
would exhibit chaos.
We will now examine the implications of such a scenario.

\section{Observations Revisited}\label{s-revisit}

By performing a more detailed numerical study, we show below that, while
Lorenz's scenario is the correct basic description of the transition, the
IC actually breaks apart before the ``cusp parameter'' value is reached.
The transition from smooth IC to chaotic attractor thus occurs over a
range of parameter values, rather than at a single critical parameter
value.
We emphasize that Lorenz did not claim that his scenario {\it did}
happen, only that there could be no simpler scenario.
Our study shows that the proposed ``simple'' scenario does not happen in
this particular family, nor should it be expected in the general transition
from smooth invariant circle to chaotic attractor in the presence of
noninvertibility.

Before beginning the description of our numerical investigation, we
recall some terminology and preliminary results.

\subsection{Noninvertible preliminaries}\label{ss-prelims}

Let $f:\Re^2 \to \Re^2$ be a smooth ($C^\infty$) map. 
Let $x$ be a period-$q$ point for $f$.
As for invertible maps, we define the {\it stable manifold} of $x$ to be
$W^s(x) = \{y \in \Re^2 : f^{kq}(y) \to x {\rm \ as \ } k \to \infty\}$.
Since the inverse map is not necessarily uniquely defined, we
define the {\it unstable manifold} of $x$ as
$W^u(x) = \{y \in \Re^2 : {\rm \ there\ exists\ a\ biinfinite\ orbit\ }
\{y_j\} \ {\rm \ with\ }
y_{i+1}=f(y_i), y_0=y {\rm \ and\ } y_{-kq}  \to x
{\rm \ as \ } k \to \infty\}$.
The use of the term ``manifold'' is an abuse of terminology for
noninvertible maps since, while the {\it local}
stable and unstable manifolds are guaranteed to be true manifolds,
the {\it global}
manifolds are not \cite{McS, Robinson, FGKM}.
(Some authors address this problem by using ``set'' instead of ``manifold''
\cite{FGKM,RMG,McN, Nien, Sander99}.)
By iteration of the local smooth
unstable manifold, the global unstable manifold can be computed as
for invertible maps.
If the manifold (or any curve) ever crosses the critical set ($J_0$, defined
below) tangent to the zero eigenvector of the map at the crossing point, its
image can have a cusp \cite{FGKM, Sander}.
More degenerate cases are also possible, but not considered here.
The global unstable manifold can also have self intersections.
Thus, it is globally a smoothly immersed submanifold, with smoothness
violated only at forward images of these special crossing points.
This is discussed further below.
The global stable manifold, due to multiple preimages, can be
disconnected, fail to be smooth, and even increase in dimension
\cite{Sander, FGKM}.
The multiple preimages lead to interesting basins
of attraction (see Section~\ref{s-discussion} and \cite{AK,MPGKC}).

A set $S \subset \Re^2$ is {\it invariant} if $f(S)=S$.
When $S$ is invariant and a topological circle we call it an
{\it invariant circle} (IC).
An invariant circle on which there exists a periodic orbit is
called a {\it frequency locked circle} or {\it circle in resonance}.
(The terminology is borrowed from
return maps of periodically forced systems.)
Typically we see a single attracting
and a single repelling periodic orbit on a frequency locked circle, but
multiple orbits are certainly possible.
A periodic orbit is also called a {\it periodic locking} or a
{\it locked solution}, although we usually reserve the use of
periodic locking to refer to a periodic
orbit {\it on} an invariant circle.

We have already indicated above that we are using the term {\it
attractor} informally to indicate whatever is suggested via computer
simulations.
More formally, we define a compact set $K$ to be an {\it attractor block}
if $f(K) \subset K^\circ$ (the interior of $K$).
Then we define $\Lambda$, the {\it attracting set} associated with
attractor block $K$ by
$\Lambda = \bigcap_{n=1}^\infty f^n(K)$ \cite{ACHM, Robinson}. 
In the cases we consider, including for $\tau$ in the range depicted in
fig.~\ref{observ}, we take the attractor block $K$ to be an annulus
containing the ``attractor'' in each picture: inner radius a small
circle around the repelling fixed point $O$, and outer radius a
(topological) circle which is big enough to contain the attractor, but
still in the right half plane.
{\it The} attracting set will always be taken to be associated with this
annulus.

Any invariant circles or periodic points in the attractor block will also
be contained in the attracting set $\Lambda$.
So will the unstable manifolds of the periodic points.
The attractors displayed by computer simulation are generally (an
approximation of) some subset of $\Lambda$.
If $\Lambda$ exhibits at least some approximate recurrence property, such as
a dense orbit (for example, a quasiperiodic IC), or a numerical orbit which
appears dense (for example, a locked circle with a high period locking),
then the attractor might be a good approximation of all of $\Lambda$.
In other cases, such as a locked circle with a low period locking,
simulation would only display the sinks.
Computation of the unstable
manifolds of the saddles (crucial to our numerical investigation)
often gives a better picture of $\Lambda$ than does a simulation.
Note that a locked circle is still an invariant circle, even though the whole
circle is not the attractor.  The circle is smooth by the unstable manifold
theorem \cite{McS} through the saddles, but as with invertible maps,
its smoothness where two
branches of the unstable manifold meet at the sinks depends on the
eigenvalues \cite{ACHM}.
Unlike invertible maps, as indicated above, the smoothness of an IC
can also be violated by the development of cusps.

\begin{figure}
\centerline{ \psfig{figure=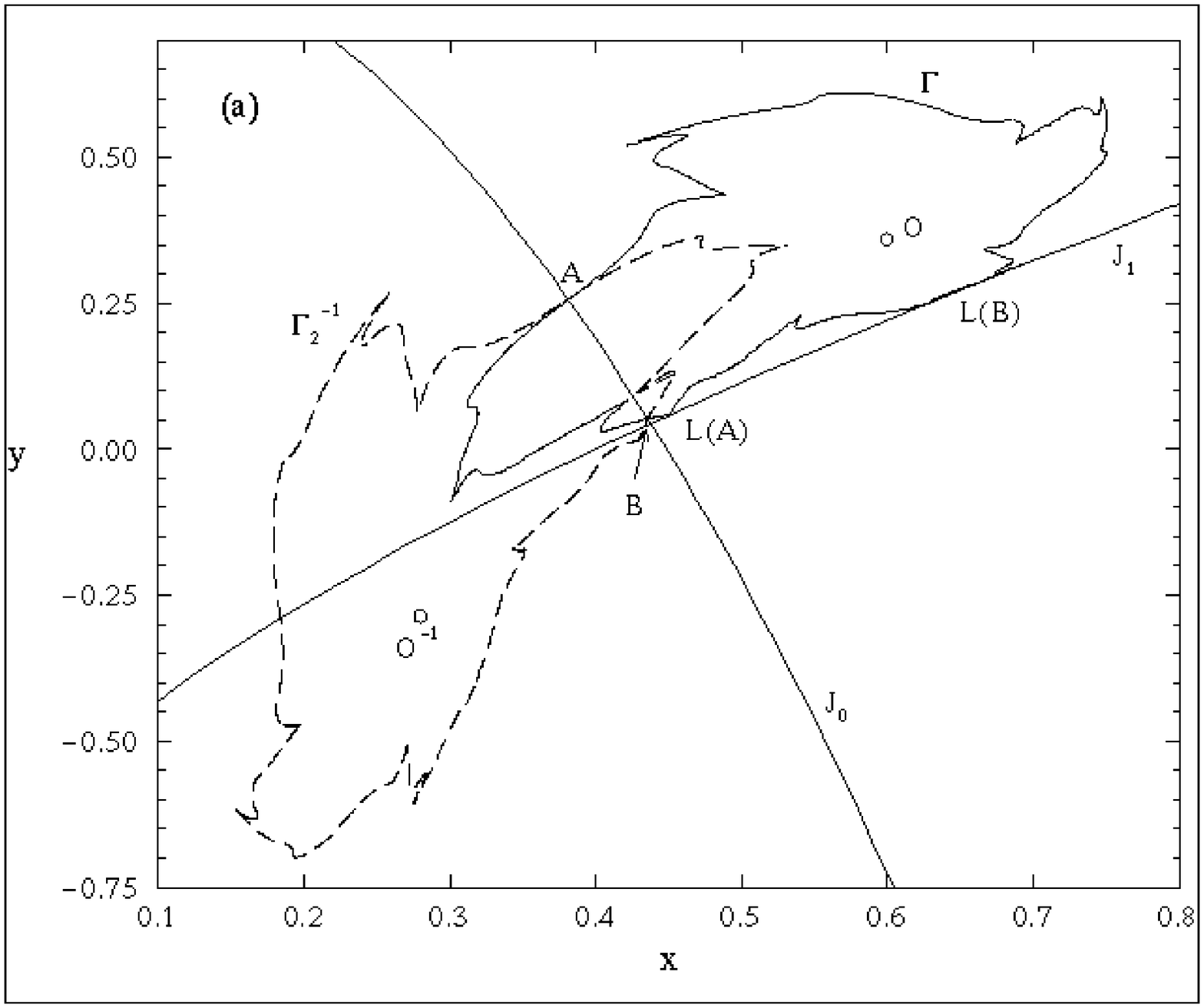,height=15pc}  
\psfig{figure=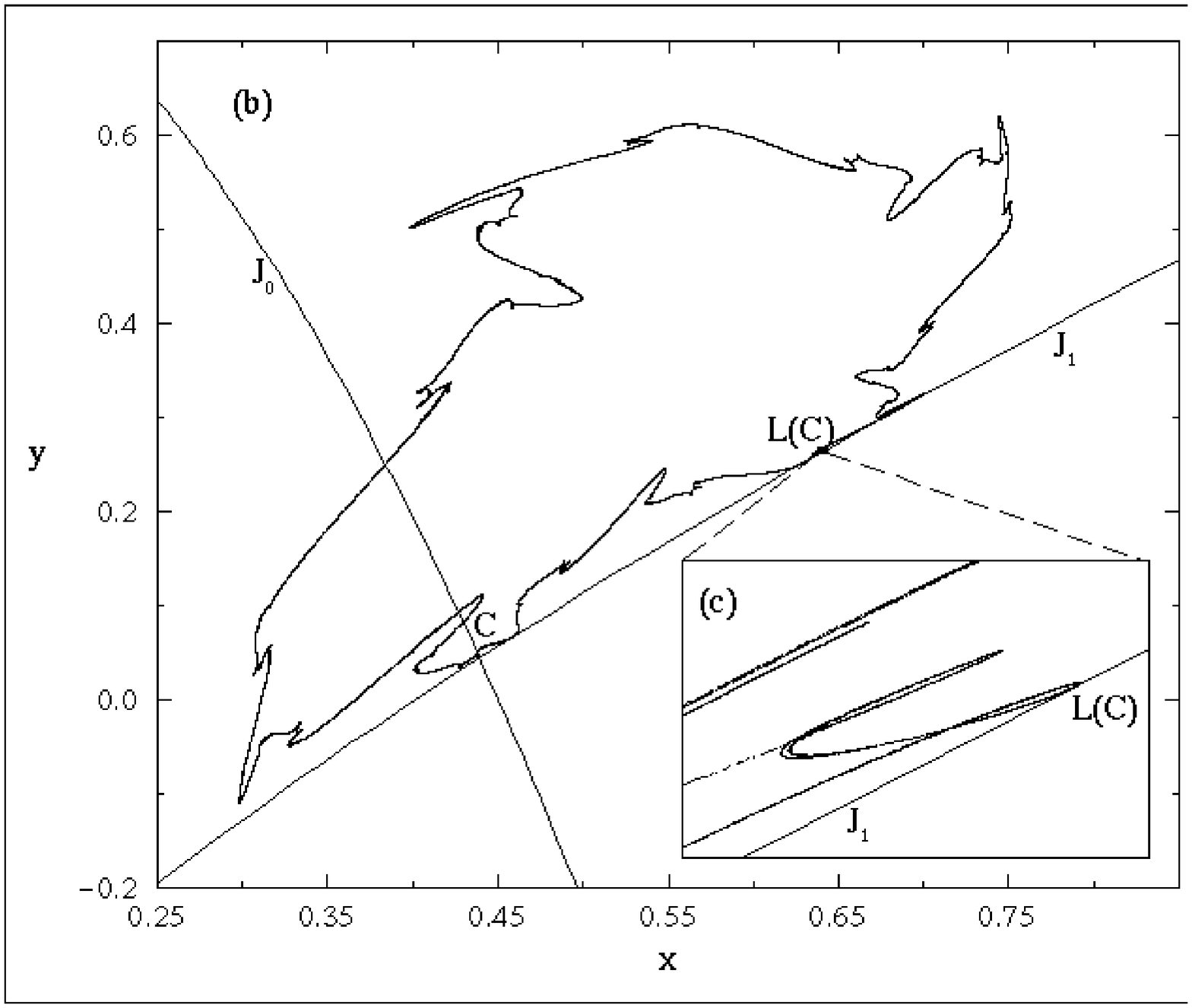,height=15pc} }
\caption{(a)
Invariant circle, $\Gamma$, one of its preimages, $\Gamma_2^{-1}$,
and the critical curves $J_0$ and $J_1$ ($\tau=1.775$).
The image of a ``curve'' crossing $J_0$ (at $A$ and $B$) is tangent to
$J_1$ (at $L(A)$ and $L(B)$, respectively).
(b) ``weak chaotic ring'' possessing loops ($\tau=1.785$); 
(c) enlargement showing the loop on the attractor at $L(C)$.
\label{gamma_inv}}
\end{figure}

We now define the {\it critical curve}
$J_0$ (termed ``curve C" in \cite{Lo2} and 
often $LC_{-1}$ for {\bf L}igne {\bf C}ritique in the recent noninvertible
literature \cite{FGKM,MiraBook3}) as the locus in phase space where the
linearization of the map becomes singular.
We call its image $J_1$ (``curve D" and $LC$, respectively).
The two sets $J_0$ and $J_1$, as well as additional images and preimages,
are known to
be a key to understanding the dynamics of a noninvertible map.
Consider, for example, Fig.~\ref{gamma_inv}a showing
an enlargement of the IC,
$\Gamma$ for $\tau=1.775$, with three important additions:
$J_0$, $J_1$, and $\Gamma_2^{-1}$,
one {\it additional} first-rank preimage of the invariant circle $\Gamma$.
The map $L_{(0.36, 1.775)}$ (abbreviated as $L$)
has either one or three inverses depending
on the phase-plane point 
in consideration (the term ``$Z_1 - Z_3$ map" has been proposed for
such maps, \cite{MPGKC}). 
The invariant circle $\Gamma$ has {\it three} first-rank
preimages: itself, the curve $\Gamma_2^{-1}$ shown in the figure,
and a third preimage $\Gamma_3^{-1}$, further away in phase space (not shown).
The geometry of the map can be visualized by first folding the left side
of the phase plane along $J_0$ and onto the right side.
$\Gamma$ and $\Gamma_2^{-1}$ should now coincide.
Next rotate roughly $90^\circ$ so that $J_0$ maps to $J_1$.
The image of $\Gamma$ lands exactly on
itself (before the folding), sending $A$ to $L(A)$ and $B$ to $L(B)$ in
the process.
Note that $L$ maps points that are in the intersection of the two regions bounded
by $\Gamma$ and $\Gamma_2^{-1}$ to the region bounded by the pieces of $\Gamma$
and $J_1$ which connect at $L(A)$ and $L(B)$.
That is, some points inside $\Gamma$ are mapped outside $\Gamma$.
Similarly, points inside $\Gamma_2^{-1}$, but not $\Gamma$ are mapped
from outside $\Gamma$ to inside $\Gamma$.
Contrast this with the property that interiors of invariant circles for
(orientation preserving) invertible maps are invariant.

Note also that $\Gamma$ is tangent to $J_1$ at $L(A)$ and $L(B)$.
The ``folding'' along $J_0$, typically creates a quadratic tangency along
$J_1$ at 
images of curves crossing $J_0$, unless the curve crossing $J_0$ does so
at an angle tangent to the zero eigenvector at that point.
This phenomenon has been extensively discussed in \cite{FGKM}.
Our explanation here is based on the understanding reached in that paper. 
As a parameter (in this case $\tau$) is varied, the map
$L_{(0.36,\tau)}$ (again abbreviated as $L$), the entire attracting set
$\Lambda$, the curves $J_0$ and $J_1$ as well as the intersection points
of $\Lambda$ with $J_0$ (e.g. $A$) and with $J_1$ (e.g. $L(A)$) vary.
%
In our case, between $\tau$ values for figs.~\ref{gamma_inv}a and b, the
invariant curve develops a ``fjord' which pushes across $J_0$, creating
two new intersection points of $\Lambda$ with $J_0$, and two new tangencies
with $J_1$ at the two image points.
We call the lower intersection point of the fjord $C$ (see
fig.~\ref{gamma_inv}b).
The points $C$ and $L(C)$ vary with $\tau$ as well.
When the tangent of the ``curve'' in $\Lambda$
at the intersection point $C$ becomes coincident with the null vector of
the Jacobian of the map at $C$, the image of $\Lambda$ ($\Lambda$ itself)
acquires a cusp touching $J_1$ at $L(C)$ \cite{FGKM}.
As the parameter value $\tau$ is further varied, the cusp becomes again
a quadratic tangency, but $\Lambda$ acquires a loop
(figs.~\ref{gamma_inv}b and c).
Iteration of this ``loop"
gives rise to an infinity of such loops
(called ``antennae" in \cite{Lo2}) on $\Lambda$, suggesting chaotic behavior.

Further, Lorenz gives a heuristic justification of a simple criterion for an
attractor (not formally defined, but assumed by us to be
a compact invariant set with some sort of recurrence) to be a chaotic 
attractor (an attractor exhibiting sensitivity to initial conditions):
the existence of two distinct points on the attracator which
map to the same image point (also on the attractor) \cite{Lo2}.
In this context, the tracing out of an object with a single orbit suggests
recurrence, and existence of a loop implies the existence of
two points which map to the same image point.
Thus, if this criterion is correct,
attractors with loops would necessarily be chaotic attractors.
%

%
%


\subsection{Numerical investigation}\label{ss-numerical}

%

\begin{figure}
\centerline{ \psfig{figure=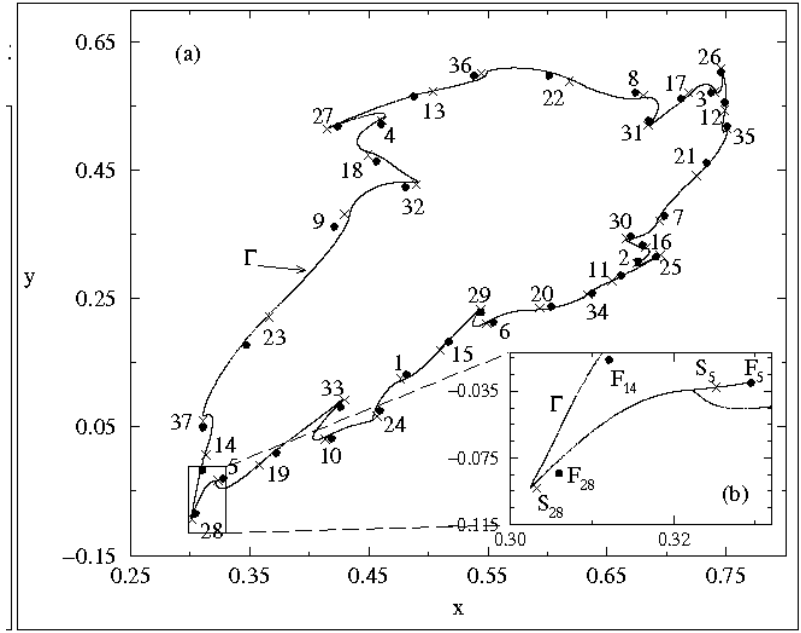,height=20pc} }
\caption{ Coexisting attractors for $\tau=1.7765$: the attracting period-$37$
foci ($\bullet$), and the IC.  The period-$37$ saddles ($\times$) are
also shown, along with both branches of the unstable manifold $S_5$.
One branch approaches the focus $F_5$, while the other branch approaches the IC.
\label{t1.7765}}
\end{figure}

\begin{figure}
\centerline{ \psfig{figure=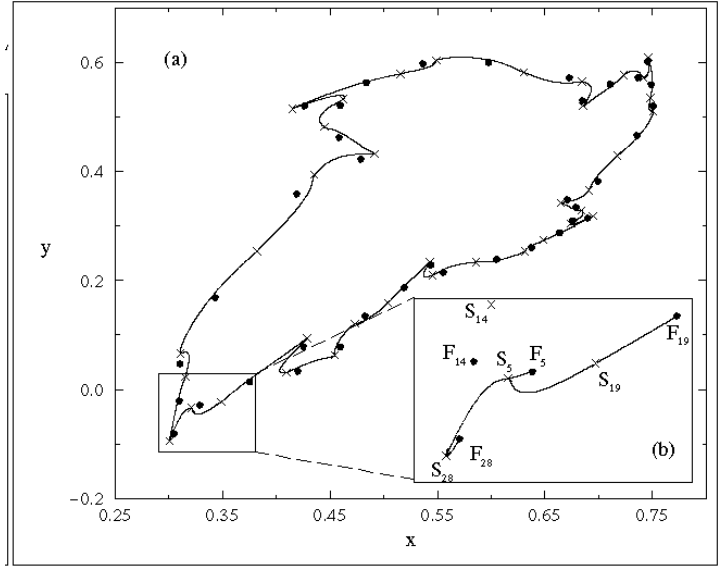,height=18pc} }
\caption{During the ``first'' manifold crossing ($\tau = 1.77688$).
a) the period-$37$ saddle and sink 
solutions have been plotted with one branch of the global unstable manifold
(the branch involved in the crossing) which emanates from each of the 37 saddles.
b) Both branches of the unstable manifold emanating from the saddle
labelled $S_{19}$
have been plotted (but {\it not} the branches emanating from $S_5$,
or any other saddle).
\label{t1.77688}}
\end{figure}

\begin{figure}
\centerline{ \psfig{figure=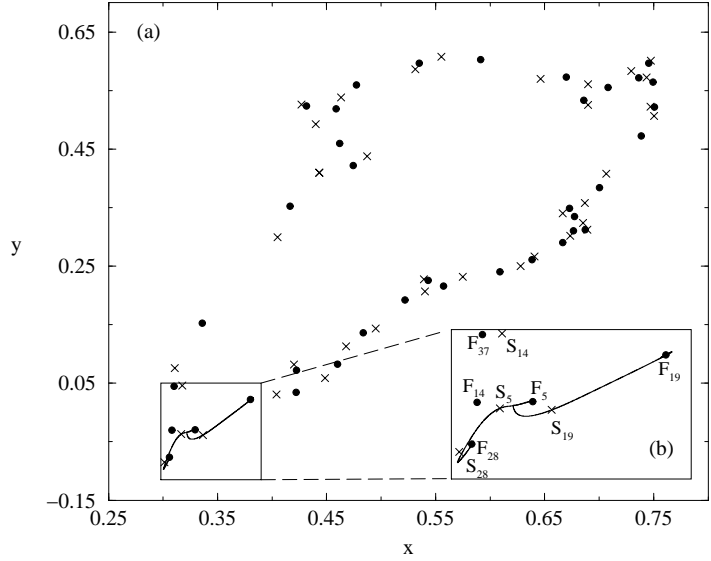,height=18pc} }
\caption{After the first manifold crossing ($\tau=1.778$): the
invariant circle as an attractor has disappeared, but a new invariant
frequency locked circle has formed as an attracting {\it set}.
Both branches of the unstable manifold of the two adjacent period-$37$ saddles
($S_{5}$ and $S_{19}$) are plotted.
The right branch of the unstable manifold of $S_5$ and the left branch of
the unstable manifold of $S_{19}$ meet at the attracting focus $F_{5}$.
\label{t1.778}}
\end{figure}

\begin{figure}
\centerline{ \psfig{figure=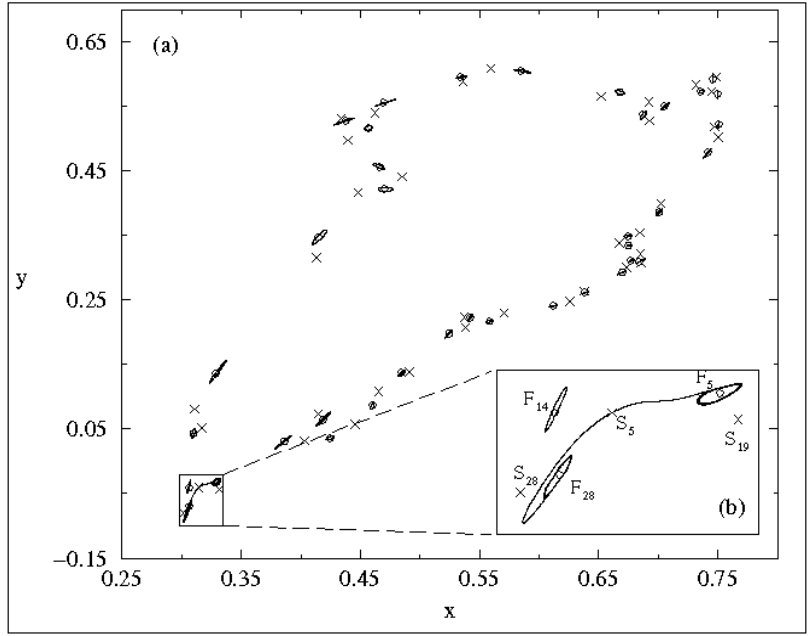,height=20pc} }
\caption{After a secondary Hopf bifurcation of the period-$37$ foci:
period-$37$ invariant circles surround unstable foci ($\circ$) at  $\tau=1.78$.
Both sides of the unstable manifold of the period-$37$ saddles --- shown here
only for $S_5$ ---
asymptote to the stable period-$37$ invariant circles.
\label{t1.78}}
\end{figure}

\begin{figure}
\centerline{ \psfig{figure=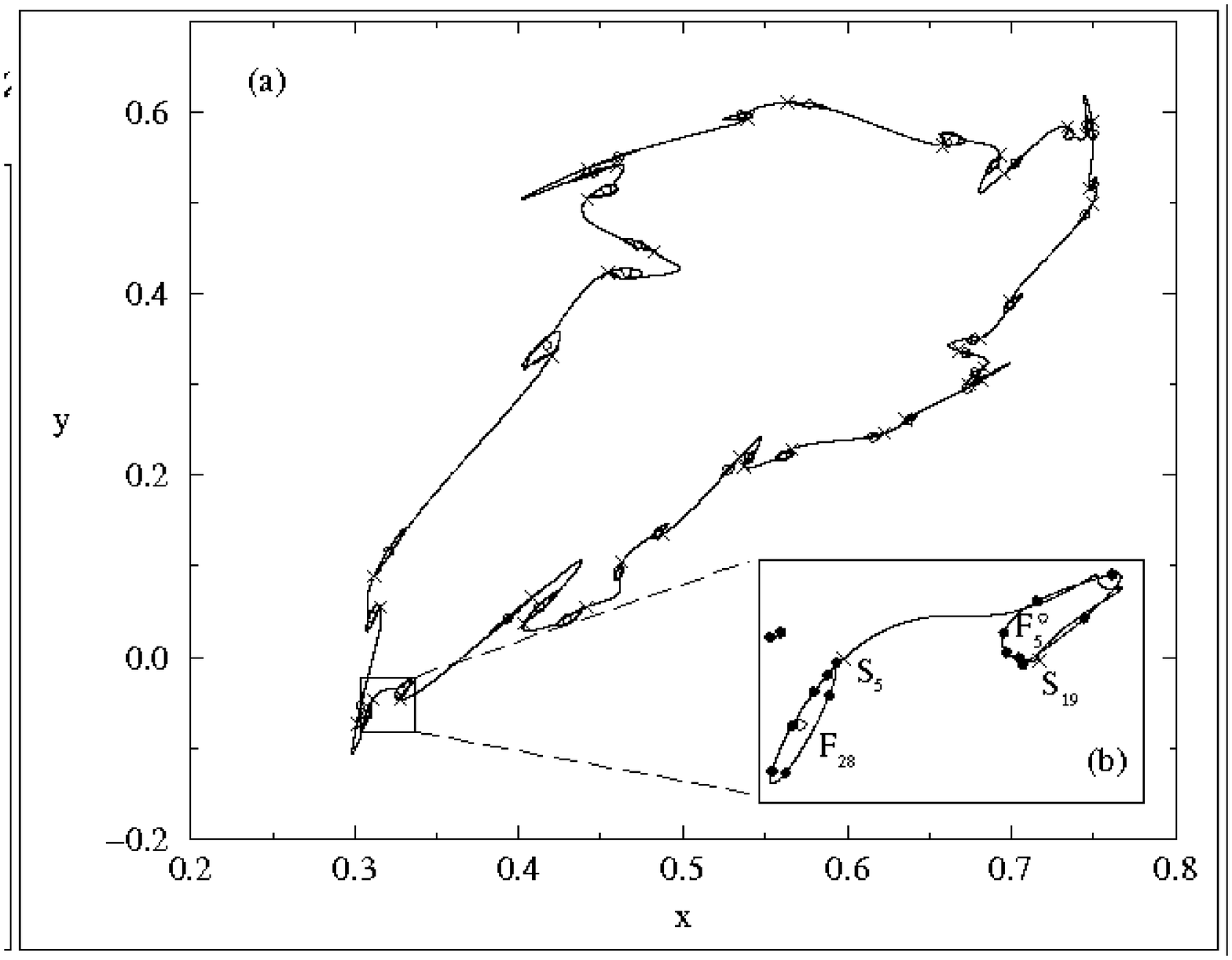,height=20pc} }
\caption{Locking on the small IC ($\tau = 1.7835$).
Both branches of the unstable manifold of $S_5$ are plotted.
Self-intersections are apparent on the right branch.
Also shown are period-$37$ saddles
($\times$), period-$37$ repelling foci ($\circ$), and, in (b),
period-$259$ sinks ($\bullet$), but not the period-$259$ saddles.
\label{t1.7835}}
\end{figure}

\begin{figure}
\centerline{ \psfig{figure=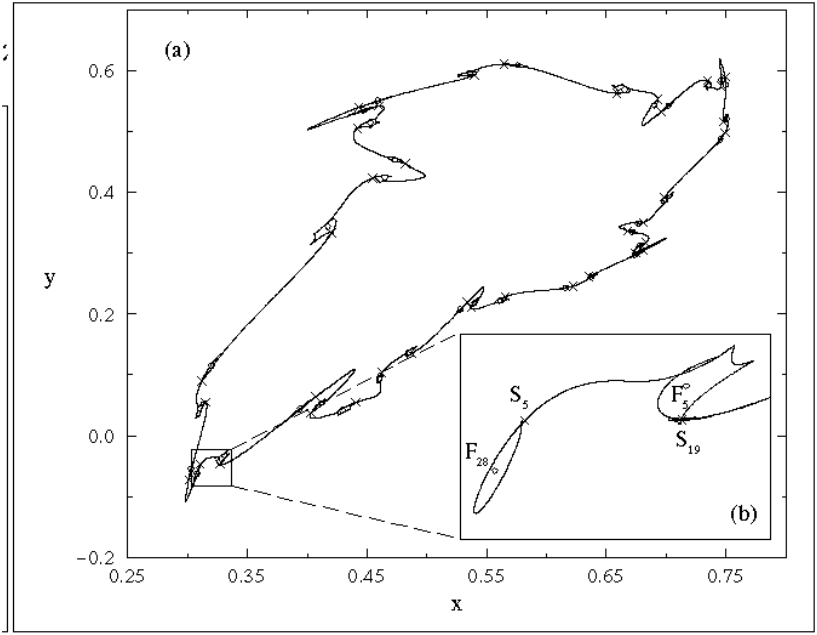,height=20pc} }
\caption{During the ``second'' manifold crossing ($\tau= 1.784$).
(a) Unstable manifolds of all $37$ saddles are plotted.
(b) Only the two branches of the unstable manifold of $S_5$ are plotted;
the right branch is the heteroclinic branch.
For clarity, neither the period-$296$ saddles nor foci locked on the small ICs
are plotted.
\label{t1.784}}
\end{figure}

\begin{figure}
\centerline{ \psfig{figure=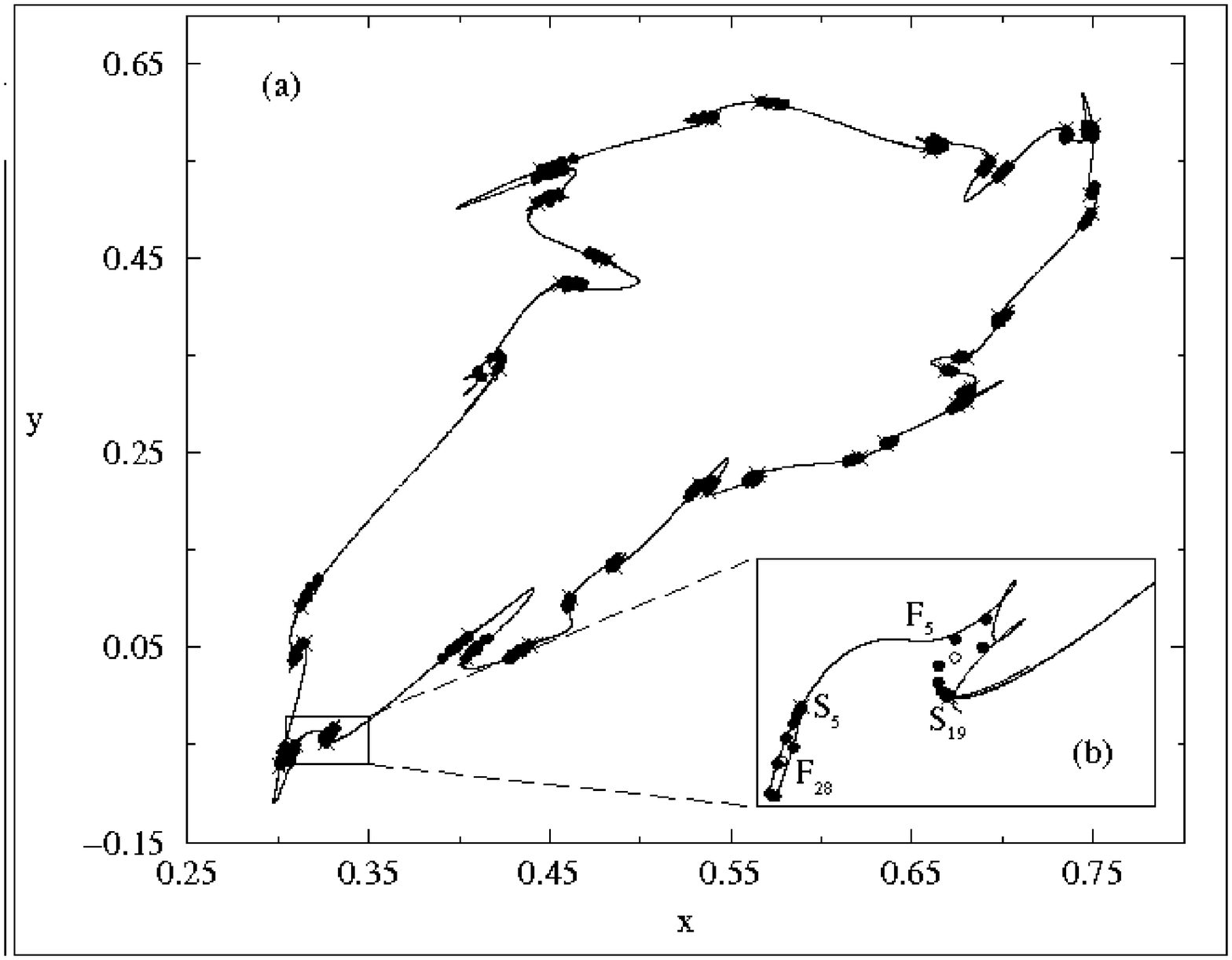,height=15pc}
\psfig{figure=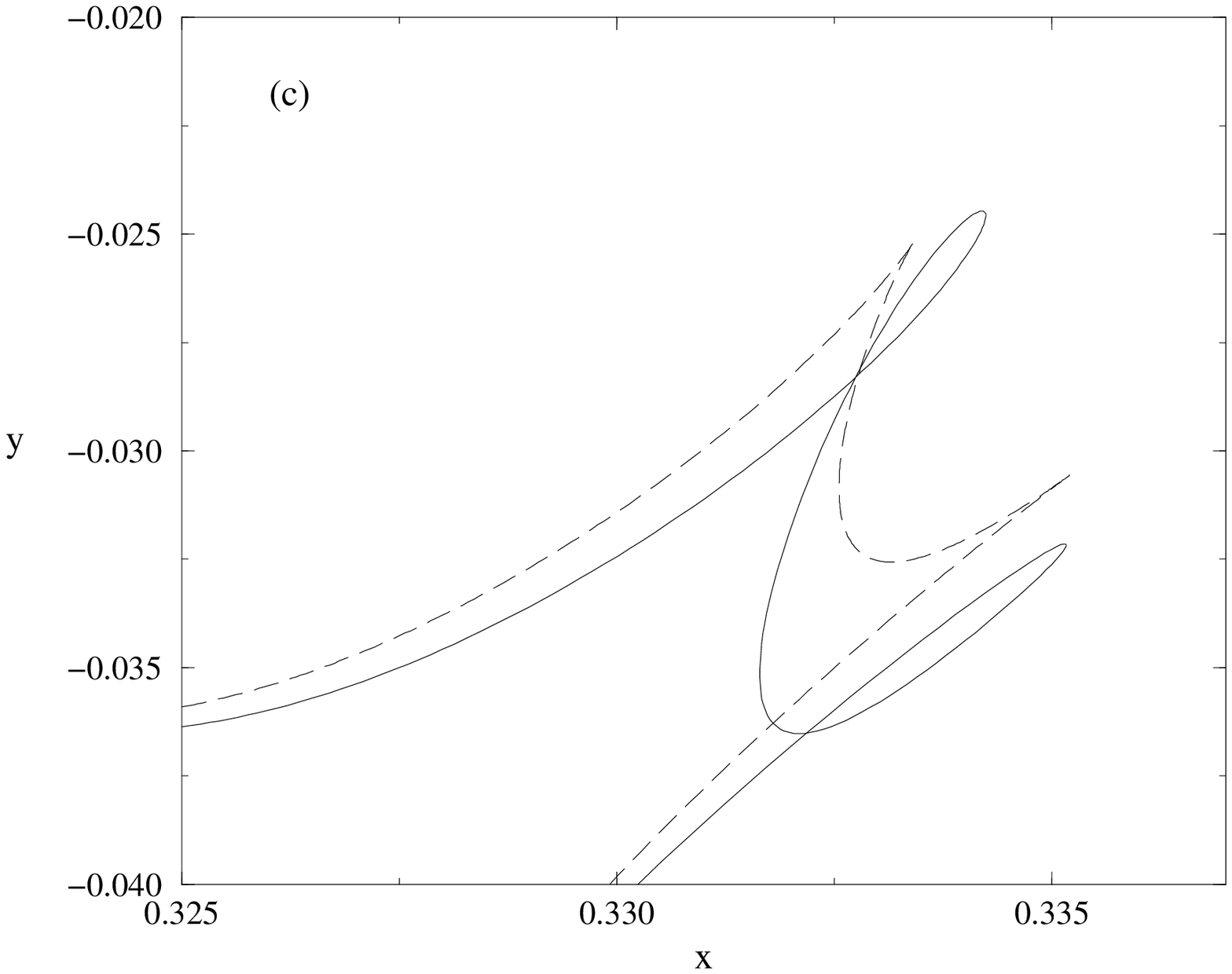,height=15pc}  }
\caption{ (a) Coexistence of period-$407$ attracting orbit with
a ``weakly chaotic ring'' ($\tau=1.785$).
Plotted are the period-$37$ saddles ($\times$), period-$37$ repelling
foci ($\circ$), period-$407$ attracting nodes ($\bullet$), and the
unstable manifolds of the period-$37$ saddles.
(b) enlargement with only the two branches of the unstable manifold
of $S_5$ plotted.
(c) Further enlargement of a portion of the saddle global unstable manifold
in (b) at $\tau = 1.785$ (solid line) with loops;
the analogous unstable manifold for the approximate ``cusp parameter'':
$\tau=1.78428$ (dashed line).
Not shown are the period-$245$ saddles and nodes which exist on the small
ICs at this cusp parameter value.
\label{t1.785}}
\end{figure}

\begin{figure}
\centerline{ \psfig{figure=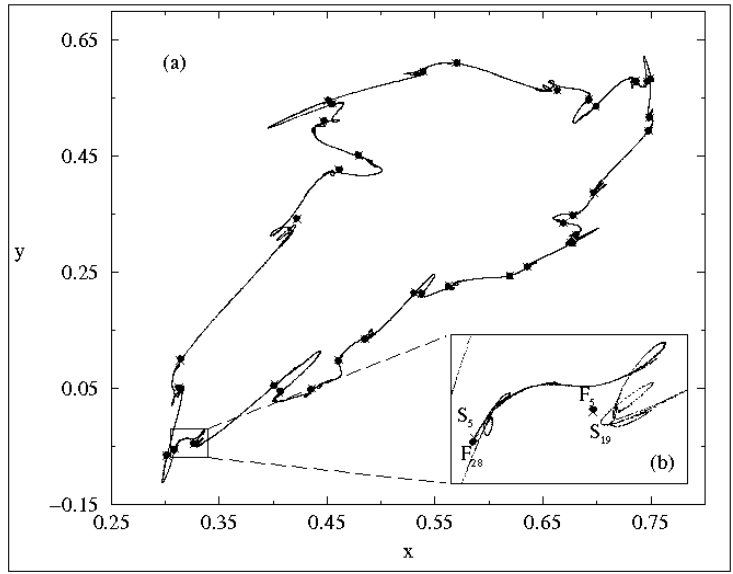,height=18pc} }
\caption{Coexistence of period-$37$ attracting orbit with
a ``weakly chaotic ring'' ($\tau=1.7866$).
(a) Plotted are the period-$37$ saddles ($\times$), period-$37$
attracting nodes ($\bullet$), and the
unstable manifolds of the period-$37$ saddles.
(b) Enlargement with only the two branches of the unstable
manifold of $S_5$ plotted.
The left branch immediately approaches the period-$37$
stable orbit ($F_{28}$);
the right branch approaches the ``big" 
attractor: it iterates off the figure to the right, follows counterclockwise
around the attractor seen in (a), and returns in the bottom of (b),
coming in just to the right of $F_{28}$ and $S_5$;
its loops are larger than in fig.~\ref{t1.785}b.
\label{t1.7866}}
\end{figure}

\begin{figure}
\centerline{ \psfig{figure=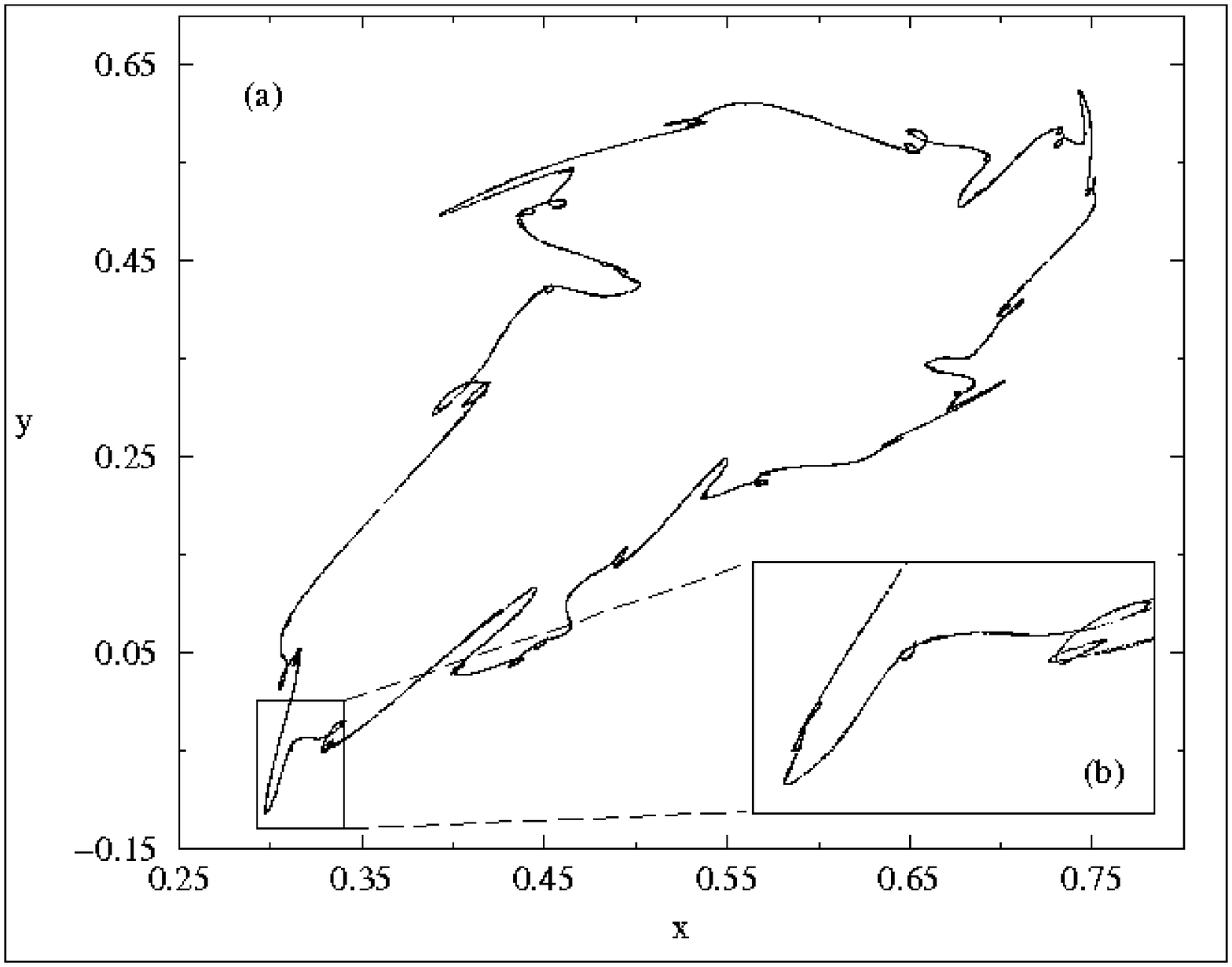,height=20pc} }
\caption{Only the chaotic attractor remains ($\tau=1.788$).
It was approximated by computing a single orbit.
The period-$37$ fixed points have disappeared in a saddle-node
bifurcation. \label{t1.788}}
\end{figure}

We set out to examine the transition from smooth IC to chaotic attractor
in finer detail in figures \ref{t1.7765} 
to \ref{t1.788}.
Before we embark on a detailed description, we should state 
that the  real picture is yet more complex, involving similar 
phenomena on an even finer scale.
For example, the same phenomena which cause the destruction of the ``large'' IC
could also destroy the ``small'' IC's which we refer to below.
Nevertheless, this sequence of figures includes the
essential features of the transition: coexistence of a periodic orbit with a
nearby IC, destruction of the IC through a global bifurcation, appearance
of loops on an unstable manifold, and the reappearance of an attractor,
this time chaotic with loops.
It might be useful to look ahead to the two-parameter bifurcation
diagram in fig.~\ref{lorbif} to put the one-parameter cut described in
the current section in context.
%
%

We start from the smooth IC at $\tau=1.775$ (fig.~\ref{gamma_inv}a).
At $\tau \approx 1.776243$ a period-$37$ saddle and node pair appears in
a saddle-node (turning point) bifurcation. 
The bifurcation occurs {\it off} of the IC, and
results in the coexistence of two 
attractors: the stable invariant circle, $\Gamma$, and an attracting 
period-$37$ sink.
Although the sinks have real eigenvalues at
the saddle-node bifurcation, the eigenvalues very quickly become complex
and they have been labelled in fig.~\ref{t1.7765}b for $\tau=1.7765$
by $F_i$ for ``focus.''
%
%
The separatrix between these two coexisting attractors 
consists of the
stable manifolds of the period-$37$ saddles.
%
%
%
%
%
Global {\it unstable} manifolds have been systematically
constructed as (37th) forward iterates of 
local unstable manifolds as for invertible maps.
In this case, we have plotted only the two branches emanating from
saddle $S_5$.
The right branch of the global unstable manifold directly approaches $F_5$.
This is clear in the blowup of fig.~\ref{t1.7765}b.
The left branch very quickly approaches the IC $\Gamma$.
Successive 37th iterates travel slowly in a clockwise direction around the IC,
eventually showing a good approximation to the whole IC.
The approach to the IC is somewhat clearer in fig.~\ref{t1.7765}b,
where the branch appears to ``join'' the IC just to the left of $S_5$.
The unstable manifold then continues almost ``on'' the IC:
past $F_{28}$, $S_{28}$, and $F_{14}$, evolving around the
IC, reentering from the right, just below $F_5$,
before completing the (approximation to the) IC just to the left of
$S_5$.

Notice also that 18 of the period-$37$ saddle-node pairs are ``inside" 
(those marked by numbers inside $\Gamma$) and the rest ``outside" 
the invariant circle; this is another feature of noninvertible maps \cite{FGKM}).
We discuss this again in Sec.~\ref{s-discussion}.

As $\tau$ increases further, the consequences of {\it global bifurcations}
are observed;
these bifurcations involve interactions between the stable
and unstable manifolds of the saddle period-$37$ points. 
The first such global bifurcation is evident in the sequence of figures
\ref{t1.7765} through \ref{t1.778}:  the left branch of the unstable manifold
of $S_{19}$ misses $S_{5}$ to the left in fig.~\ref{t1.7765} and approaches
the IC (the curve coming into the inset from the right and below $F_5$ is
practically coincident with this unstable manifold branch); it makes the
``saddle connection'' with $S_5$ in fig.~\ref{t1.77688}; it misses $S_{5}$
to the right in fig.~\ref{t1.778} and approaches $F_5$.

The global crossing actually ``starts" with a heteroclinic tangency
(between a branch of a stable manifold
of a period-$37$ saddle and and branch of an
unstable manifold of another saddle
on the same orbit) which destroys the IC
at $\tau \approx 1.776878$.
This tangency is also called a crisis \cite{OttBook},
since as $\tau$ increases
toward the tangency parameter
value, the IC (approached by the unstable manifold) and its basin
boundary (the stable manifold) approach each other, colliding at the
tangency parameter value.
The attracting set during the crossing appears --- as for invertible
maps --- to be a complicated topological object which at least includes
the closure of the unstable manifolds of the period-$37$ saddles.
Numerically, most orbits appear to approach the period-$37$ sink
(fig.~\ref{t1.77688}).
Thirty-seventh iterates of each of the $37$ branches of the unstable
manifold plotted in fig.~\ref{t1.77688}a
evolve clockwise toward the ``next'' saddle, suggesting that we are at least
close to a heteroclinic crossing of stable and unstable manifolds.
Although the stable manifolds of the saddles cannot
be realistically computed because of the exponential explosion 
(``arborescence") of the number of period-$37$ preimages,
the crossing of the stable manifold is evident in the inset 
(fig.~\ref{t1.77688}b).
The left branch of the unstable manifold of $S_{19}$
heads toward saddle $S_5$, and although it is not clear from the figure,
it begins to alternate back and forth, with folds on the right
approaching $F_5$, and folds on the left approaching $S_{28}$.
The portion approaching $S_{28}$ folds again, with one side heading toward
$F_{28}$ (visible) and the other side heading toward $S_{14}$
(only the start is visible).
Further iterates (not shown) would result in this one branch of the unstable
manifold of $S_{19}$
eventually evolving all the way around the
attracting set suggested in fig.~\ref{t1.77688}a, repeatedly
visiting every saddle and every node.
As with invertible maps, this allows for long chaotic transients, if one
could appropriately choose initial conditions on the unstable manifold.
Most orbits, however, eventually are attracted to the period-$37$ sink.

We note that a transverse crossing of stable and unstable manifolds in the
noninvertible setting does not imply transverse crossings along the rest
of the orbit.
For this reason, homoclinic or heteroclinic crossings do not
necessarily imply the
existence of a ``horseshoe'' or the chaotic transients referred to above.
However, Sander \cite{Sander} shows that the failure of a transverse
crossing to provide a horseshoe is a codimension-two phenomenon, so we
assume that our transverse crossings in our one-parameter study
avoid such a problematic point, and behave as in the invertible case.
Thus we expect the parameter value at
the tangency initiating this heteroclinic crossing
to be an upper bound on $\tau_{CD}$, marking the first appearance
of chaotic dynamics.
The invariant set with the chaotic dynamics
does not appear to be an attractor, however, so
the chaos is numerically evident only as transient behavior.

After the crossing, which ends in a heteroclinic tangency at
$\tau \approx 1.776881$,
both branches of the unstable manifolds of the period-$37$ saddles
approach the period-$37$ sink.
The only attractor that remains is the stable period-$37$ sink, shown along
with the period-$37$ saddles in fig.~\ref{t1.778} for $\tau=1.1778$.
The attracting {\it set}, however, is a ``new'' invariant circle, with a
period-$37$
locking, formed by connecting the saddles and nodes with the unstable manifolds
of the saddles.
%
%
This invariant circle is seen in the blowup of fig.~\ref{t1.778}b.
%
%
It is smooth except at the sinks,
where two different branches connnect
and the smoothness depends on the eigenvalues \cite{ACHM}.
Here the IC is only $C^0$ since the sinks are foci, with complex eigenvalues.
There are apparently no cusps on the unstable manifolds at this parameter value.

As $\tau$ is increased, a secondary Hopf bifurcation of the period-$37$ foci at
$\tau=1.779444$ renders them unstable. This can be seen in
fig.~\ref{t1.78} for 
$\tau=1.78$, where the unstable foci, are surrounded
by stable period-$37$ invariant circles.  
It is important (and barely visible in the figure) that the 
unstable manifold {\it intersects} the small invariant circles
(and itself!) repeatedly. 
This is a reminder of the noninvertibility of the map, and our
acknowledged abuse of the term ``manifold."
See \cite{FGKM} for a discussion
of this type of unstable manifold self-intersection.

A number of high-period lockings on the period-$37$ invariant circles
appear for very narrow ranges of the parameter. 
For example, at $\tau=1.7835$, a period-$259$ saddle-node pair 
(formed after a period-7 locking of the period-$37$ invariant circles)
is visible (fig.~\ref{t1.7835});
assuming the period-$7$ locking is {\it on} the small ICs, these small
ICs are now frequency locked circles, composed of the chain of the 
unstable manifolds of the period-$259$ saddles (not shown). 
The period-$37$ unstable manifolds asymptotically approach 
the small ICs.
The right branch clearly intersects itself and the IC it
approaches (fig.~\ref{t1.7835}b).
Most individual orbits eventually approach the attracting period-$259$
sinks.

A second heteroclinic
manifold crossing involving the ``other" side of the 
global unstable manifold of the
period-$37$ saddle ``starts'' at $\tau \approx 1.78373$.
The bifurcation is apparent in figures \ref{t1.7835} through \ref{t1.785}:
the right branch of the unstable manifold
of $S_{5}$ swings around $F_5$, and misses $S_{19}$ to the left in
fig.~\ref{t1.7835}, approaching the small IC around $F_5$;
it makes the ``saddle connection'' with $S_{19}$ in fig.~\ref{t1.784};
it misses $S_{19}$ to the right in fig.~\ref{t1.785},
with successive iterates following essentially ``on'' the attractor with loops.

As for the earlier manifold crossing, understanding the saddle connection at
$\tau = 1.784$ requires some explanation.
At this parameter value
(fig.~\ref{t1.784}), a period-$8$ resonance on the small invariant circles
results in the period-$296$ points apparently being the only attractors.
The left branch of the unstable manifold of $S_5$ very quickly approaches
the small IC surrounding $F_{28}$ and continues counterclockwise around it;
it is not a homoclinic loop, although it appears close to being one.
The right branch heads for $S_{19}$ after looping  clockwise about $3/4$ of
the way around $F_5$ (including the two high curvature ``corners'');
the global crossing is apparent since part of the manifold iterates to the
left of $S_{19}$ around the IC (we iterated only half way around the IC),
and part of the manifold iterates to the right of $S_{19}$ toward the right edge
of the inset;
this part to the right of $S_{19}$ continues the
``long transient'' behavior.
It is important to notice that although both branches of the unstable
manifold of $S_{19}$ are being followed by parts of the right branch
of the unstable manifold of $S_5$, neither branch of the unstable
manifold of $S_{19}$ has been directly plotted in this figure. 

The slope of the unstable manifold at the point of intersection 
with $J_0$ continues to change as $\tau$ is increased.
At $\tau = \tau_{cusp}\approx 1.78428$
the slope becomes equal to the slope of the eigenvector corresponding to the
zero eigenvalue, and {\it cusp points} appear on the ``manifold" 
(fig.~\ref{t1.785}c).
Beyond this ``cusp parameter'' value, the unstable manifold contains loops.
Since $\tau_{cusp}$ occurs while the unstable manifold is still
intersecting the stable manifold, the picture for $\tau$
just beyond $\tau_{cusp}$, but not beyond the heteroclinic tangency
which marks the end of the stable/unstable manifold crossing, is more
complicated than might be expected in other families.
In our case, the newly developed loops are damped out as they eventually
approach the small period-$37$ ICs (possibly with locked solutions on the
small ICs).
Of course, the loops could follow a long transient before being attracted
to the ICs.
We do not have a figure corresponding to such a parameter value.

With further increase of $\tau$, we pass the ``other end'' of the manifold 
crossing (beyond the final heteroclinic tangency) and an apparently
chaotic attractor appears.
Thus, this final heteroclinic tangency parameter value appears to be
an upper bound on $\tau_{chaos}$, marking 
the first appearance of a chaotic attractor.
At $\tau = 1.785$ the resulting attractor is seen in fig.~\ref{t1.785} 
to coexist with a stable period-$407$ (born at a 
period-$11$ locking on the period-$37$ small invariant circles):
the left branch of the unstable manifold of $S_5$
approaches the nearby small IC with the period-$407$ locking;
the right branch approaches the large amplitude ``messy" attractor that
self-intersects 
and exhibits the loops it has inherited from the global
unstable manifold.
Such an attractor has been called a ``weakly chaotic ring''
(\cite{FGKM,MiraBook3}). 
Note that Lorenz plotted this weakly chaotic ring attractor by iterating
a single orbit starting near the fixed point (fig.~\ref{observ}c);
we plotted the same attractor by following a branch of the unstable
manifold of the period-$37$ saddle orbit (fig.~\ref{t1.785}).
Tracking the periodic saddles and their manifolds
also allowed us to see the coexisting
small IC with the  period-$407$ attractor which was not apparent in
fig.~\ref{observ}c.

The small amplitude invariant circles now begin to shrink,
disappearing at $\tau \approx 1.78626168$ in
a second secondary Hopf bifurcation;
the period-$37$ foci from return from sources to sinks;
these continue to coexist  with the large 
attractor that has become very ``messy" with many more points of 
self-intersection (fig.~\ref{t1.7866}b for $\tau=1.7866$). 
The stable period-$37$ foci approach the period-$37$ saddles in anticipation of 
the following period-$37$ saddle-node bifurcation, which happens at
$\tau \approx 1.7866262$.


Finally, at $\tau=1.788$ the saddle-node bifurcation of the period-$37$ 
solutions has occurred, and the pair has disappeared.
Note that (as with invertible maps) the pairings have switched.  For
example, $S_5$ disappears with $F_{28}$, while $F_5$ disappears with
$S_{19}$.
Only the large chaotic attractor is left, seen in fig.~\ref{t1.788}. 

\par{\bf Recap.}
The transition from smooth IC attractor
($\tau=1.75$, figs.~\ref{observ}b, \ref{gamma_inv}a) to
noninvertible ``messy" chaotic attractor ($\tau=1.788$, fig.~\ref{t1.788})
is not a result of a smooth change 
of a smooth, large amplitude invariant circle which develops cusps
at a critical parameter value.
Instead, the transition occurs via a sequence of bifurcations, notably
involving the interaction of the IC with nearby
periodic orbits and the unstable manifolds of the saddles.
Specifically, the IC as an attractor is destroyed in 
a global manifold crossing that
leaves an attracting periodic
orbit as the only attractor.
The global unstable 
manifolds of the saddle period-$37$ points are the invariant objects that
first develop cusps (at $\tau_{cusp} \approx 1.78428$,
fig.~\ref{t1.785}c)
and later loops (beyond $\tau_{cusp}$).
%
When the cusps and loops first appear, they ``damp out'' as they approach
the attracting periodic orbit (or small IC).
When the large amplitude attractor reappears at after the second global
manifold crossing, it is approached by the branch of the unstable manifold
which now has loops.
These loops are inherited by the large amplitude attractor which is now
apparently chaotic.

\vfill \eject

\section{The larger picture.}\label{s-larger}

\begin{figure}
\centerline{ \psfig{figure=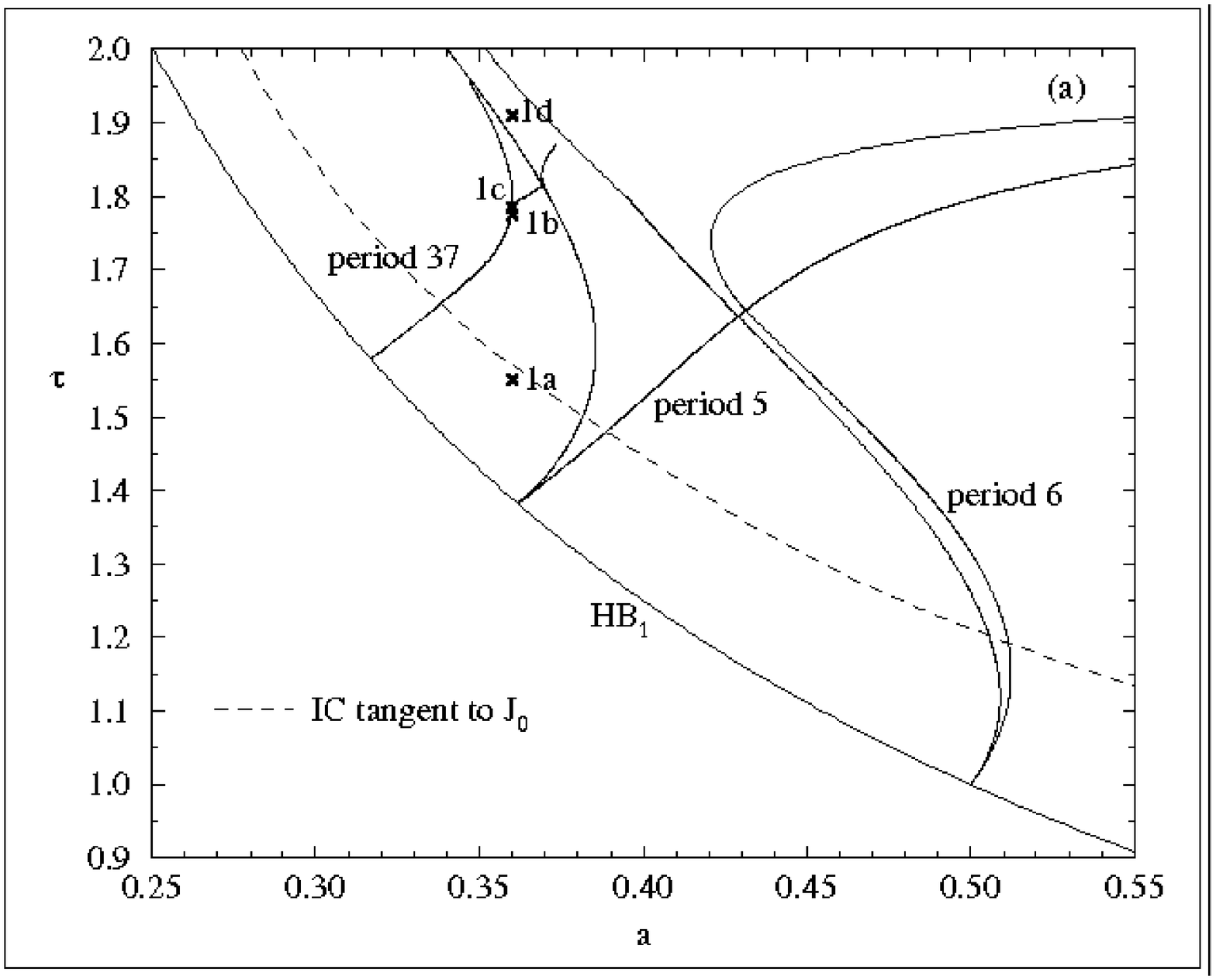,height=16pc}
\psfig{figure=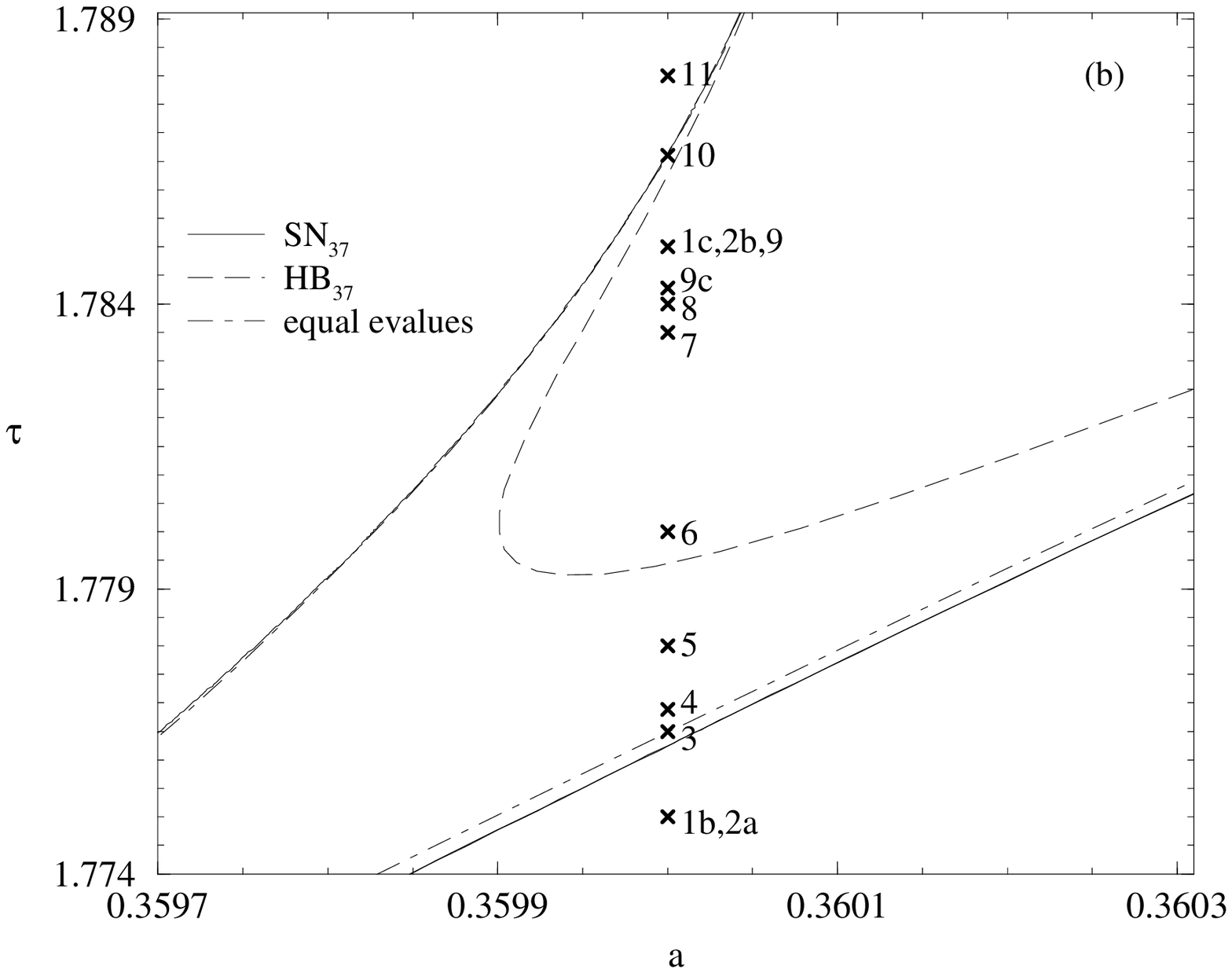,height=16pc} }
\caption{(a) The start of a two-parameter bifurcation diagram for the map
$L_{(a,\tau)}$, containing the fixed-point Hopf bifurcation curve HB$_1$,
and three pairs of saddle-node curves emanating from HB$_1$ corresponding to
period-$5$, period-$6$, and period-$37$ resonance horns.
The $\times$'s are labelled by the phase portraits in fig.~\ref{observ}.
(b) An enlargement of the period-$37$ resonance horn where it is traversed
in our one-parameter cut.
The two period-$37$ saddle-node curves (SN$_{37}$) mark the upper and lower
boundaries of the resonance horn.
Two equal eigenvalue curves and a secondary Hopf curve (HB$_{37}$)
subdivide the horn.
The parameter values marked with $\times$'s are labelled by our
corresponding figure number(s).
 \label{lorbif}}
\end{figure}

While a rich dynamical picture occurs along our 
one-parameter family with respect to $\tau$, it is
clear that it only constitutes a part of a larger, actually
two-parameter scenario.
This is consistent with existing two-parameter descriptions of Hopf bifurcations
for invertible maps \cite{Arnold, ACHM, Chenciner, Pthesis, P, PFK}.
The primary features in the invertible case are resonance horns
--- inside which the corresponding maps have a periodic orbit of a
certain period ---
emanating from the primary Hopf bifurcation curve.
%
Since the fixed point at a  Hopf bifurcation
has eigenvalues on the unit circle, the IC close to its birth lies
{\it far away} in phase space from $J_0$, where the map has a zero eigenvalue.
Thus, we expect the two-parameter bifurcation diagram for the noninvertible
Hopf bifurcation to strongly resemble the analogous bifurcation diagram for
the invertible case --- at least near the Hopf bifurcation curve.

We have embedded our one-parameter family ($a=0.36$ in eq.~(\ref{map})) in
the two-parameter family $L_{(a, \tau)}$.
Fig.~\ref{lorbif}a shows part of the ``big picture"
in the $\tau$-$a$ plane:
the period-one Hopf bifurcation curve together with period-$37$, 
period-$5$, and period-$6$ resonance horns emanating from it.
The bifurcation points we encountered along the one-parameter cut
``continue'' to bifurcation curves in the two-parameter setting.
The curves displayed in this figure are only a tiny sampling
of the dozens of bifurcation curves we have computed,
which are, in turn, only
a tiny sampling of the bifurcation curves that
can be computed for this family.
(In particular, there are saddle-node curves in addition to those emanating
from the primary Hopf bifurcation curve.
See \cite{MP, PK} for
descriptions of such intraresonance horn features for invertible maps.)

As we move higher up in the horns
--- above the ``IC tangent to $J_0$'' curve in fig.~\ref{lorbif}a ---
the circle grows and
starts to interact with $J_0$ and thus becomes ``truly" noninvertible.
It is above this curve where the corresponding maps can
exhibit noninvertible behavior, such as
mapping some of the interior of an IC to its exterior and vice versa (recall
fig.~\ref{gamma_inv}a).
More ``interesting'' behavior, such as coexisting attractors or the breakup
of the IC, is suggested when the resonance horns start to overlap.
In our one-parameter cut this overlap happens at least by the first
period-$37$ sadde-node bifurcation at $\tau \approx 1.776243$, although
the three horns we have displayed do not overlap until higher $\tau$
values.

The enlargement in  fig.~\ref{lorbif}b shows three of the multitude of
bifurcation curves which
subdivide the period-$37$ horn into ``inequivalent'' regions: two
``equal eigenvalue'' curves (just inside the horn),
marking the transition between real and complex
eigenvalues for a period-$37$ orbit, and a secondary Hopf curve
(HB$_{37}$), marking the birth of the ``period-$37$ small ICs.''
Somewhat unsatisfying in fig.~\ref{lorbif}b is the fact that
none of the displayed bifurcation curves is unique to noninvertible maps.
We remark that we did compute an ``eigenvalue zero'' curve
(for the existence of a period-$37$ point with an eigenvalue of zero),
but we did not display it in fig.~\ref{lorbif}b because it had no
direct bearing on the bifurcation
sequence we described (but see such a curve in a period-$6$ horn in
fig.~\ref{pe6}).
However, this curve certainly has an indirect bearing on the bifurcations,
because it is related to the crossing of $J_0$ by the
attracting set (which includes any ICs, periodic orbits, and their
unstable manifolds).
If the attracting set never interacted with $J_0$, the
bifurcations could be produced with invertible maps.

Other bifurcation curves, which our current computational tools did not
handle easily, would separate the other
$\times$'s in fig.~\ref{lorbif}b from each other.
For example, the $\times$ labelled \ref{t1.77688}
is in a thin heteroclinic region; heteroclinic tangency curves separate it
from \ref{t1.7765} below and \ref{t1.778} above.
Similarly, the $\times$'s labelled \ref{t1.784} and \ref{t1.785}c are
in a thin heteroclinic region, separating them from \ref{t1.7835}
below and \ref{t1.785} above.
A ``cusp curve'' (unique to {\it noninvertible} resonance horns) runs through
\ref{t1.785}c.
The $\times$ labelled \ref{t1.7835} is in a thin ``period-$259$
secondary resonance horn," which emanates from the secondary Hopf curve
HB$_{37}$ and corresponds to a locking on the small period-$37$ ICs.

A more complete knowledge of the internal structure
of the ``noninvertible'' resonance horns in the two-parameter setting
is clearly the key to extending
our current understanding of the 
transition from a smooth invariant circle to a chaotic attractor.

\section{Discussion} \label{s-discussion}

\subsection{Typical one-parameter cuts}\label{ss-typical}

A natural question which arises following the detailed numerical
investigation of Sec.~\ref{s-revisit} is how ``typical'' are the
features of that specific one-parameter family.
We take the liberty in the next few paragraphs of using the intuition we
have  gained via our numerical explorations, including the partial
two-parameter bifurcation diagram of the previous section, to provide
at least a partial answer to this question.
Even though somewhat speculative, we feel that these observations are
useful for placing the one-parameter cut in a broader context.

A key observation accompanying our understanding of the IC breakup is that
an IC with a cusp {\it and} an irrational rotation number would have the cusp
iterate forward to a dense set of image points, which would also be cusps.
Even if such an object is topologically possible, it appears to us to be
highly unlikely.
Because of this, we expect the invariant circle to break up {\it before}
it reaches a parameter value where it develops a cusp.
Thus, the IC must interact with some other object in the attracting set.
In our experience, the other object in this ``pre-cusp'' breakup is a nearby
periodic orbit and its accompanying stable and unstable manifolds.
The  related manifold crossings appear to be the ``natural'' way to first
destroy an existing IC and to later reconstitute ``itself'' as a chaotic
attractor.

This observation is consistent with the sequence of transitions we described
in  Sec.~\ref{s-revisit}.  Extending this one-parameter transition interval
to the two-parameter setting, we imagine a transition band, roughly parallel
to the
`IC tangent to $J_0$' curve, and passing between parameter values labelled
`1b' and `1d' in fig.\ref{lorbif}a.
Below the lower boundary of the transition band a smooth invariant circle
exists (possibly locked, possibly coexisting with periodic orbit(s)).
Above the upper boundary of the transition band, the corresponding maps are
chaotic.
Loosely, the transition band is a ``thick version'' of the
``cusp curve'' suggested by Lorenz.

Further supporting this transition band scenario is the prominent ``opening''
of the resonance horns in fig.~\ref{lorbif}a.
This would be even more prominent if more resonance horns were plotted:
most of the resonance horns appear to open up and overlap near the top
of the upper boundary of the transition band.
This is suggestive of, or at least consistent with, chaos.

The existence of the transition band also suggests why the period-$37$
orbit might
have been the specific
orbit that was involved in our one-parameter scenario: that was the horn
which was passed through as the transition to chaos band was crossed.
Vertical parameter cuts for values of $a$ other than $a=0.36$ could pass
through a different periodic horn during this narrow transitional
parameter interval.
We note that Lorenz specifically chose $a=0.36$ because it appeared to
``avoid periodic windows.''
In the two-parameter context, this can be
interpreted as avoiding low-period resonance horns.
Thus, we might expect a ``typical'' one-parameter cut to often
interact with a periodic orbit of a period less than $37$.

Alternative viable bifurcation sequences are suggested by the blowup of the
period-$37$ horn in fig.~\ref{lorbif}b.
For example, moving the cut to the left of $a=0.3599$ could avoid the
secondary Hopf bifurcation curve HB$_{37}$, and thus avoid the small
invariant circles that appeared in the $a=0.36$ cut.
Other differences in varying $a$ would be more obvious if we were able
to compute the region(s) of the period-$37$ horn corresponding to the
global bifurcations, and curves corresponding to the development of
cusps.  (We expect to be able to do this in the future.)
We know the $\times$ labelled `$9c$' in fig.~\ref{lorbif}b is both
on a cusp curve (in the parameter space) and in a heteroclinic region.
Moving the vertical cut either left or right would probably prevent these two
from coinciding.
That is, the fact that the cusp parameter value occured during the manifold
crossing in the $a=0.36$ cut is probably not typical of most
one-parameter cuts.
If the cusp parameter value occured first (with respect to $\tau$), then
the phase portrait for the first appearance of loops would be simpler:
the loops would simply approach the small attracting ICs.
(If we encountered the cusp parameter outside the secondary Hopf curve,
it would be still
simpler: the loops would damp out at the attracting period-$37$ nodes.
See fig.~\ref{damp} for comparison.)
On the other hand, if the cusp parameter value occured after the
heteroclinic region had been crossed, we claim (with no
proof) that there must be another periodic orbit interacting with the IC
as well.  Perhaps {\it its} unstable manifold will be the one to
develop a cusp and {\it later} undergo a global manifold crossing to
reconstitute the chaotic attractor with loops.

\subsection{Related noninvertible features}\label{ss-related}

The relatively high period of the locking encountered along
the one-parameter cut above makes the elucidation of the phase and
parameter space structures cumbersome.
To illustrate some of these phenomena more succinctly we perform
computations in the neighborhood of a simpler (period-$6$) resonance
horn.

\begin{figure}
\centerline{ \psfig{figure=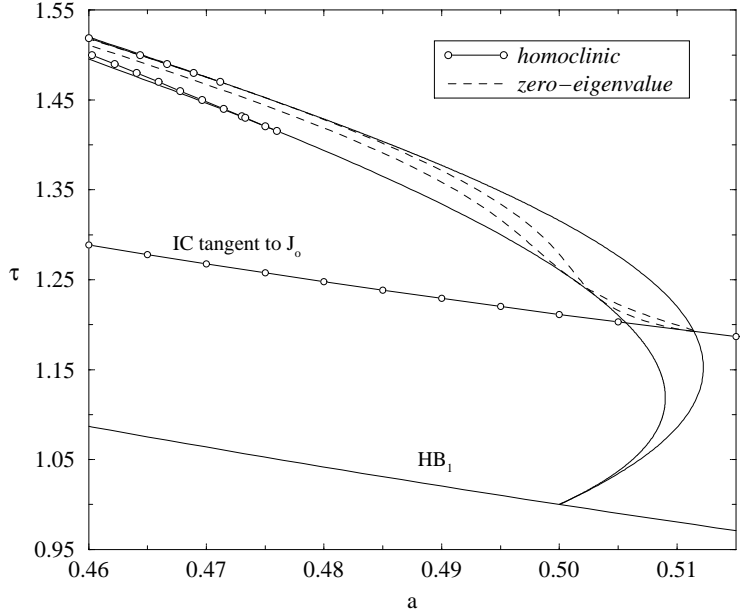,height=20pc} }
\caption{ Period-$6$ resonance horn and related bifurcation curves.
\label{pe6}}
\end{figure}

Fig.~\ref{pe6} is a blow-up of this period-$6$ resonance horn, with
two additional types of bifurcation curves plotted inside: 
(a) two global bifurcation curves (approximating a thin heteroclinic manifold
crossing region associated 
with the stable and unstable manifolds of the period-$6$ saddle points)
and (b) the curve on which the periodic points acquire an eigenvalue
zero (i.e. one periodic point lies on $J_0$ and one on $J_1$)
\cite{GFK, GAK}.
This latter curve provides an indication of the interaction of the
dynamics with noninvertibility.
An additional curve (`IC tangent to $J_0$'), extending beyond
the resonance horn, is plotted.
This is the locus, in parameter space,
where the invariant circle starts interacting with noninvertibility,
by first touching the critical curve $J_0$.
Note that the `eigenvalue zero' curve is tangent to the `IC tangent to $J_0$'
curve.

\begin{figure}
\centerline{ \psfig{figure=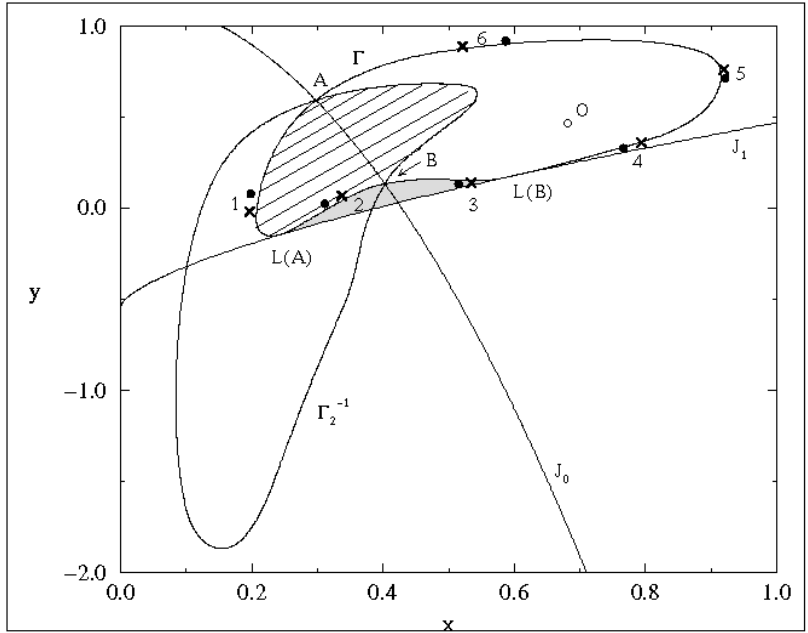,height=20pc} }
\caption{Inside the horn: Period-$6$ saddle-node pairs inside and outside the
invariant circle ($a=0.465$, $\tau=1.472$); the hatched region inside $\Gamma$
is 
mapped to the grey region outside. \label{in_out6}}
\end{figure}

In fig.~\ref{in_out6},
we show a phase portrait for a parameter value inside the period-$6$
horn and above the `IC tangent to $J_0$' curve ( $(a,\tau) = (0.465, 1.472)$ ).
The invariant circle {\it crosses} $J_0$ transversely at
two points $A$ and $B$;
it is easy to see that it will 
then become quadratically tangent to $J_1$ at the
image points $L(A)$ and $L(B)$  (see \cite{FGKM,MiraBook3}).
Also illustrated in the figure is the
relative location
of the attracting invariant circle and one more of its three
first-rank preimages.
Clearly, one of the period-$6$ saddle-node pairs is located {\it inside} the 
attracting invariant circle while the remaining five are located {\it outside}.
(Note that iterates of points to the left of $J_0$ switch sides of the
IC;  thus, there must be an even number of periodic points to the left
of $J_0$.)
We have also illustrated the portion of the interior of the
invariant circle that
is mapped ``outside" after one iteration of the map, as well as the portion of
phase space outside the invariant circle that gets mapped inside.

\begin{figure}
\centerline{ \psfig{figure=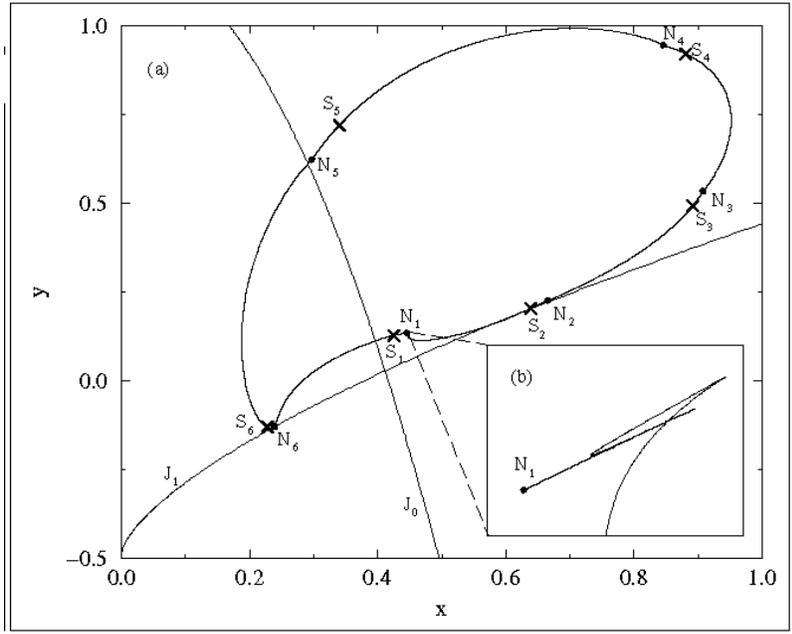,height=20pc} }
\caption{Simple damping out of a cusp at an attracting period-$6$
node ($(a,\tau)=(0.4838, 1.409)$ ).
(b) The upper right cusp maps (under $L^6_{(0.4838, 1.409)}$) to what
looks like the end of a line segment emanating from $N_1$.  The whole unstable
``manifold'' has
effectively collapsed onto this line.  Successive (sixth) iterates of the
cusp approach $N_1$ along the line. }\label{damp}
\end{figure}

Our next observation is for another
parameter value inside the period-$6$ horn:
$(a,\tau)=(0.4838,1.409)$.
In fig.~\ref{damp}
we have an example of an unstable manifold of a period-$6$ saddle with a
cusp which damps out at the period-$6$ attracting orbit.
Recall that in our one-parameter cut, the cusp on the period-$37$ unstable
manifold developed in the middle of a heteroclinic crossing.
Iterates of the cusp could therefore follow a long transient before
eventually being damped out as it approached a period-$37$ sink.
Fig.~\ref{damp} is a simpler version of the damping out of cusps
on an unstable manifold in the approach to a sink.

\begin{figure}
\centerline{ \psfig{figure=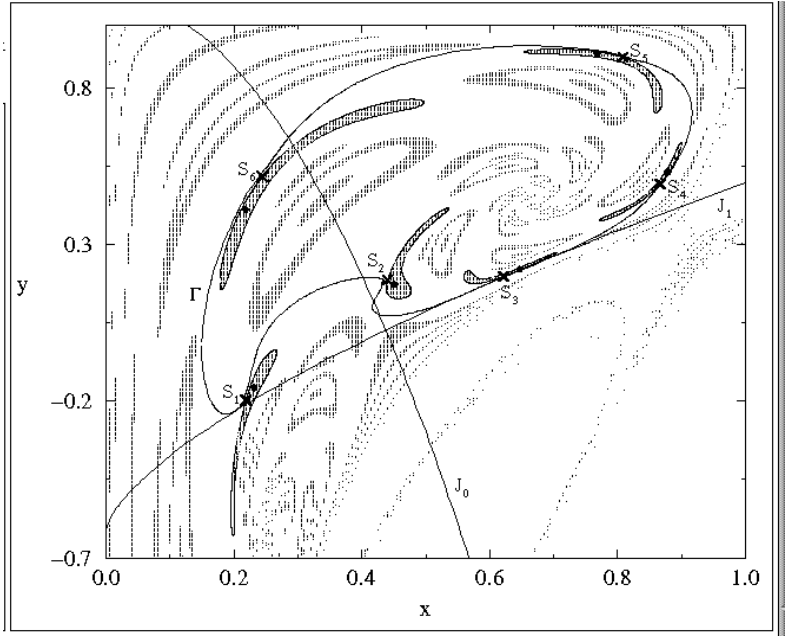,height=20pc} }
\caption{ Basin of attraction of the period-$6$ foci marked by filled circles 
(grey region) and of the coexisting invariant circle (the white complement) 
at $(a,\tau)=(0.45,1.563)$.  The computed ``immediate" stable manifolds of the 
period-$6$ saddles (marked by crosses) are also plotted. \label{pe6_bas}}
\end{figure}

Finally, we have made an attempt to illustrate the complexity of the
basins of attraction
of the coexisting attractors: the stable period-$6$ points and
the stable invariant circle.
While for the purposes of this paper we will not dwell on the details of this
complicated structure (for a discussion of the basic building blocks of such
a picture see again \cite{FGKM} as well as \cite{AK,FAK,FAKGY,GFK,GAK, Billings1}), 
we emphasize a few key elements:
\begin{itemize}
\item[(a)] the basins of attraction are not simply connected;
\item[(b)] the stable ``manifolds'' consist of a large number
(possibly, even an infinity) 
of branches (resulting from the multiple backward in time trajectories).
Only the computed ``immediate" local manifold of the period-$6$ saddles is shown 
in fig.~\ref{pe6_bas}; the images and preimages of them form the boundaries of 
the grey disconnected ``islands" that belong to the basin of attraction of the
period-$6$ foci.
\item[(c)] notice that the stable set intersects $J_0$. The images and preimages
of the intersection points have been called ``points of alternance" \cite{MiraBook2}
on parts of the stable manifolds: these are points (on $J_0$ and its preimages) 
where two preimages of a branch of the stable set come together; the ``direction 
of relative movement" of nearby points switches at such points of alternance.
\end{itemize}

\subsection{Computational issues} \label{ss-computation}

As a final comment on the current state of our work on noninvertible systems,
we note that computational tools for the accurate approximation and
visualization of objects such as stable and unstable manifolds
are crucial in both understanding
what is possible and in analyzing specific examples.
The possibility of multiple inverses
can cause an exponential
explosion in the number of distinct segments of an invariant curve or of
stable manifolds of periodic points.
The same is true for the critical set: all the images and preimages of the
critical curve $J_0$.
In three dimensions the computational issues are even further complicated:
the $J_0$ curve becomes a surface, and its iterates and
preiterates must be constructed.
Standard one-dimensional continuation tools
have to be extended (to simplicial continuation, or through software like
PISCES \cite{Wicklin}) to effectively compute, triangulate and visualize 
these surfaces and their intersections.
Constructing scientific computing, database management and visualization tools
for these objects in two and three dimensions is an important task,
and to our knowledge, an open area of research.
Development of these tools would facilitate the analysis of
many physical systems that are inherently noninvertible: 
discrete-time adaptation,
control and decision problems, ranging from Model
Reference Adaptive Control systems to signal processing applications
-in addition to more traditional noninvertible models arising in population
dynamics and economics.
These tools will also be applicable in the case of multiple {\it forward}
trajectory systems (such as those arising in {\it implicit} integrators, 
or in model predictive control), as well as, more generally, in {\it relations} 
\cite{RMG}.

\section{Conclusions}\label{s-conclusions}

The breakup of invariant circles in invertible maps of the plane possess a 
relatively well-studied dynamic phenomenology;
homoclinic tangles of the invariant manifolds
of nearby periodic saddle points, which arise from lockings on or off
the invariant
circle, play an important role in the two-parameter picture.
Invariant circles for noninvertible maps of the plane possess, in some
sense, an even ``richer" phenomenology.
Initially, close to the Hopf bifurcation 
giving rise to them, they are {\it locally} invertible, and the Arnold
horn scenario holds.
Large amplitude noninvertible invariant circles, however, may interact
with the critical curves of the map.
This can (and does) lead to new transitions and bifurcations.
We demonstrated that the one-parameter transition discovered in \cite{Lo2}
actually  constitutes a part of the interior structure of the Arnold horns
in the noninvertible two-parameter setting;
we placed it in the appropriate context by pointing out the role
of nearby periodic points and their stable and unstable manifolds.
The picture we have presented here is far from complete,
however we hope that our work is a step toward a more complete
understanding of the transition from smooth invariant circle to
chaotic attractor in the presence of noninvertibility.

\vskip .2in
\noindent
{\bf Acknowledgments.} Besides being motivated by the original work of
Professor E. Lorenz, 
we acknowledge extensive discussions with him
when this paper was first being developed in the early and mid 1990's.
He declined, however, being a co-author of
the paper, something that we believe he deserved, and we note this here.
We also acknowledge discussions with Prof. R. Rico-Martinez of the Instituto
Tecnologico de Celaya (Mexico) and Prof. C. Mira of the
Complexe Scientifique de Rangueil, Toulouse, France.
We also acknowledge the support of
the National Science Foundation (IGK and BBP: Grant No. DMS-9973926),
and the Swiss Office of Energy (CEF).
Any opinions, findings, and conclusions or recommendations expressed in
this paper are those of the authors and do not necessarily reflect the views
of the funding agencies.
\vfill \eject


\begin{thebibliography}{xxxx}
%
\bibitem{MiraBook4} Abraham R H, Gardini L and Mira C, ``Chaos in
Discrete Dynamical Systems," Springer-Verlag (TELOS) New York (1997).  
%
\bibitem{AK} Adomaitis R A and  Kevrekidis I G, ``Noninvertibility
and the structure of basins of attraction in a model adaptive control system,"
{\it J. Non Linear Sci.}, {\bf 1}, 95-105, (1991).
%
\bibitem{Arnold} Arnold V I, {\it Geometric Methods in the Theory of Ordinary
Differential Equations}, Springer-Verlag, New York, (1983).
%
\bibitem{ACHM} Aronson D G, Chory M A, Hall G R, and McGehee R P, 
``Bifurcations from an invariant circle for two-parameter families of 
maps of the plane: A computer assisted study," {\it Commun. Math. Phys.}, 
{\bf 83}, 303-354, (1982).
%
\bibitem{Billings1} Billings L and Curry J H, ``On noninvertible
maps of the plane: Eruptions," {\it Chaos} {\bf 6} 108, (1996).
%
\bibitem{Chenciner} Chenciner A, ``Bifurcation de points fixes elliptiques.
II. Orbites p\'eriodiques et ensembles de Cantor invariants," 
{\it Inventiones Mathematicae}, {\bf 80}, 81-106, (1985).
%
\bibitem{ColletEckmann} Collet P, and Eckmann J-P , 
``Iterated maps of the interval as dynamical systems," Birkh\"{a}user,
(1980).
%
\bibitem{Feigenbaum} Feigenbaum M, ``Quantitative universality for
a class of nonlinear transformations," {\it J. Stat. Phys.}, {\bf 19}(1),
25-52, (1978). 
%
\bibitem{FAK} Frouzakis C E, Adomaitis R A, and Kevrekidis I G,
``On the dynamics and global stability characteristics of an adaptive control
system," {\it Comp. Chem. Eng.} {\bf 20} Suppl. B pp.1029-1034, (1996).
%
\bibitem{FAKGY} Frouzakis C E, Adomaitis R A, Kevrekidis I G, Golden
M P, Ydstie B E,  ``The structure of basin boundaries in a
simple adaptive control system,'' in {\it Chaotic Dynamics: Theory and
Practice}, T. Bountis, ed., 195-210, Plenum Press, NY, (1992).
%
\bibitem{FGKM} Frouzakis C E, Gardini L, Kevrekidis I G, 
Millerioux G and Mira C, ``On some properties of invariant sets of 
two-dimensional noninvertible maps," {\it Int. J. Bif. Chaos}, {\bf 7}(6), 
(1997).
%
\bibitem{GAK} Gicquel N, Anderson J S, Kevrekidis I G,
``Noninvertibility and resonance in discrete-time neural networks for time-series
processing,''
{\it Phys. Lett. A} {\bf 238}(1) (1998).
%
\bibitem{GFK} Gicquel N, Frouzakis C E, Kevrekidis I G,
``Invariant circles of plane endomorphisms: a computer-assisted study,''
CESA'96, Computational Engineering in Systems Applications, Lille,
July 9-12, (1996).
%
\bibitem{MiraBook1} Gumowski I and Mira C, ``Recurrences and discrete
dynamic systems," Springer Verlag, Singapore, (1980).
%
\bibitem{Lo1} Lorenz E N,
``Deterministic nonperiodic flow," {\it J. Atmos. Sci.}, {\bf 20},130, (1963).
%
\bibitem{Lo2} Lorenz E N,
``Computational Chaos - A prelude to computational instability,"
{\it Phys. D}, {\bf 35}, 299-317, (1989).
%
\bibitem{RMG} McGehee R P, ``Attractors for closed relations
on compact Hausdorff spaces,'' {\it Indiana U. Math. J.}, {\bf 41}(4), 
1165-1209, (1992).
%
\bibitem{MiraBook2} Mira C, ``Chaotic dynamics," World Scientific,
Singapore, (1987).
%
\bibitem{McN} McGehee R P and Nien C-H,
``The images of the singular sets for two-dimensional quadratic maps," 
preprint.
%
\bibitem{McS} McGehee R P and Sander E, ``A new proof of the stable
manifold theorem," Z angew Math Phys {\bf 47}, 497-513, (1996).
%
\bibitem{Nien} Nien C-H ``Analyticity of the center-unstable manifold," 
preprint, U. of Minnesota, (1996).
%
\bibitem{MP}McGehee R P and Peckham B B, ``Arnold flames and resonance
surface folds,''   Geometry Center Research Report GCG84 and
{\it International Journal of Bifurcations and Chaos} {\bf 6}, No. 2 (1996)
315-336.
%
\bibitem{MPGKC} Mira C, Fournier-Prunaret D, Gardini L, Kawakami H,
Cathala J-C, ``Basin bifurcations of two-dimensional nonvinvertible
maps. Fractalization of basins,"  {\it Int. J. of Bif. and Chaos},
{\bf 4}(2), 343-381, (1994).
%
\bibitem{MiraBook3} Mira C, Gardini L, Barugola A, and Cathala J-C ,
``Chaotic Dynamics in Two-Dimensional Noninvertible Maps," World Scientific,
Singapore, (1996).
%
%
\bibitem{OttBook} Ott E, ``Chaos in Dynamical Systems,'' Cambridge
University Press, New York, 1993.
%
\bibitem{Pthesis}Peckham B B, ``The closing of resonance horns for periodically
forced oscillators,'' Thesis, University of Minnesota, 1988.
%
\bibitem{P}Peckham, B B,``The Necessity of the Hopf Bifurcation for
Periodically Forced oscillators with closed resonance regions,''
{\it Nonlinearity} {\bf 3} 261-280, 1990.
%
\bibitem{PFK}Peckham B B, Frouzakis C, and Kevrekidis I G,
``Bananas and Banana Splits: A parametric degeneracy in the Hopf
bifurcation for maps,'' {\it SIAM J. Math. Anal.} {\bf 26}, No. 1, 
190-217, 1995.
%
\bibitem{PK}Peckham B B and Kevrekidis I G,
``Lighting Arnold flames: resonance in doubly forced periodic
oscillators,'' {\it Nonlinearity} {\bf 15}, 
405-428, 2002.
%
\bibitem{Robinson} Robinson C, ``Dynamical Systems: Stability, Symbolic
Dynamics, and Chaos,'' CRC Press Inc., Boca Raton, Florida, 1995.
%
\bibitem{RAK}Rico-Martinez R, Adomaitis R A, and Kevrekidis I G,
``Noninvertibility in neural networks, ''
    {\it Computers and Chemical Engineering} 24, 2417-2433(2000).
%
\bibitem{Sander99} Sander E, ``Hyperbolic sets for noninvertible
maps and relations," {\it Disc. Cont. Dyn. Systems},
{\bf 5}(2):339-358, (1999). 
%
\bibitem{Sander} Sander E, ``Homoclinic tangles for noninvertible
maps,'' Nonlinear Analysis {\bf 41} 259-276 (2000).
%
\bibitem{Wicklin} Wicklin F J, ``PISCES: a platform for implicit
surfaces and curves and the exploration of singularities," manual available
on-line at the URL http://www.geom.umn.edu/locate/pisces, (1995).

\end{thebibliography}
\end{document}